\documentclass[11pt,reqno]{amsart}

\usepackage{amssymb}

\usepackage[english]{babel}

\usepackage[dvipsnames,svgnames,x11names,hyperref]{xcolor}
\usepackage[hyphens]{url}
\usepackage[colorlinks=true,linkcolor=Maroon,citecolor=blue,urlcolor=blue,hypertexnames=false,linktocpage]{hyperref}
\usepackage{bookmark}
\usepackage{amsmath,thmtools,mathtools}
\mathtoolsset{showonlyrefs=true}
\usepackage{bm}
\usepackage{mathtools}

\usepackage{fancyhdr}
\usepackage{esint}
\usepackage{enumerate}

\usepackage[scr]{rsfso}
\usepackage{mathtools}
\usepackage{pictexwd,dcpic}
\usepackage{graphicx}

\usepackage[labelfont=bf, font=small]{caption}

\newcounter{mgncount}

\declaretheorem[name=Theorem,numberwithin=section]{thm}

\declaretheorem[name=Lemma,sibling=thm]{lemma}
\declaretheorem[name=Proposition,sibling=thm]{prop}

\declaretheorem[name=Definition,style=definition,sibling=thm]{defn}
\declaretheorem[name=Corollary,sibling=thm]{cor}

\declaretheorem[name=Remark,style=remark,sibling=thm]{bem}

\numberwithin{equation}{section}

\newcommand{\ga}{\gamma}

\newcommand{\cE}{\mathcal{E}}

\newcommand{\cZ}{\mathcal{Z}}

\newcommand{\ds}{\displaystyle}

\newcommand{\<}{\left<}
\renewcommand{\>}{\right>}

\DeclareMathOperator{\dive}{div}

\DeclareMathOperator{\tr}{tr}

\DeclareMathOperator*{\osc}{\text{osc}}

\newcommand{\vp}{\color{black}}

\newcommand{\dv}{\text{\normalfont div}}

\newcommand{\R}{{\mathbb R}}
\newcommand{\M}{{\mathcal M}}
\newcommand{\dg}{d_h}
\newcommand{\cn}{\mathcal N}

\numberwithin{equation}{section}

\begin{document}
	
	\title[ROF on manifolds]{Lipschitz regularity \\ for manifold-constrained ROF elliptic systems}
	\author{Esther Cabezas-Rivas, Salvador Moll, Vicent Pallardó-Julià}
	
	\address{E. Cabezas-Rivas: Departament de Matem\`atiques,
		Universitat de Val\`encia, Av. Vicent Andrés Estellés 19, 46100 Burjassot, Spain.
		{\tt esther.cabezas-rivas@uv.es  }}
	
	\address{S. Moll: Departament d'An\`{a}lisi Matem\`atica,
		Universitat de Val\`encia, Av. Vicent Andrés Estellés 19, 46100 Burjassot, Spain.
		{\tt j.salvador.moll@uv.es }}
	\address{V. Pallardó--Julià: Institut für Analysis und Scientific Computing, TU Wien,
Wiedner Hauptstraße 8-10,
1040 Wien, Österreich.
{\tt vicent.pallardo.julia@tuwien.ac.at}}

\date{\today. The authors have been partially supported by project PID2022-136589NB-I00  funded by
	MCIN/AEI/10.13039/501100011033 and by ERDF A way of making Europe. E. Cabezas-Rivas is partially supported by the
	project CIAICO/2023/035 funded by Conselleria d’Educació, Cultura, Universitats i Ocupació, and by  the network RED2022-134077-T funded by MCIN/AEI/10.13039/501100011033. S. Moll also acknowledges partial support by the network RED2022-134784-T funded by MCIN/AEI/10.13039/501100011033. V. Pallardó--Julià has also been supported by the Agència Valenciana de la Innovació, programme ``Promoció del talent'' (Innodocto, Ref. INNTA3/2022/16).}	
	
	\begin{abstract}
We study a generalization of the manifold-valued Rudin-Osher-Fatemi (ROF) mo\-del, which
involves an initial datum $f$ mapping from a curved compact surface with smooth boundary to a complete,
connected and smooth $n$-dimensional Riemannian manifold. We prove the existence
and uniqueness of minimizers under curvature restrictions on the target and topological ones on the range of $f$. We obtain
a series of regularity results on the associated PDE system of a relaxed functional with
Neumann boundary condition. We apply these results to the ROF model to obtain Lipschitz
regularity of minimizers without further requirements on the convexity of the boundary. Additionally, we provide  variants of the regularity statement of independent interest: for 1-dimensional domains (related to signal denoising), local Lipschitz regularity (meaningful for image processing) and Lipschitz regularity for a version of the Mosolov problem coming from fluid mechanics.
	\end{abstract}
	\maketitle

	\section{Introduction and main results}

\subsection{A regularity problem within the framework of image denoising}

Let $(\Sigma, g)$, ($\mathcal N, h$) be respectively a compact surface, possibly with  smooth  boundary $\partial \Sigma \neq \emptyset$, and  a  complete connected smooth $n$-dimensional Riemannian manifold. The latter  will be treated as an isometrically embedded submanifold $\mathcal N \hookrightarrow \R^N$ in some Euclidean space,  by Nash theorem \cite{Nash1956}. We further consider the space $C^{k, \alpha}(\Sigma; \mathcal N)$ of maps $u: \Sigma \rightarrow \cn$ whose $k$th derivatives are H\"older continuous with exponent $\alpha \in (0, 1]$.	Usually we write $\overline \Sigma = \Sigma \cup \partial \Sigma$ to stress that a claim is valid up to the boundary.

Given $f \in L^2(\Sigma;\cn)$ and  $\lambda > 0$, we study the regularity of minimizers  of the  energy functional
\begin{equation}\label{eq:rof-functional-0}
\mathcal E(u) := \int_\Sigma |d u| \,d\mu_g + \frac{\lambda}{2}\int_\Sigma \dg^2(u, f)d\mu_g,
\end{equation}
where $\dg$ denotes the geodesic distance on $\mathcal N$, and $d\mu_g$ is the volume element corresponding to $g$. Here, by considering $\cn$ as embedded in $\R^N$, the energy density  reads in local coordinates as
\[|d u|^2 = \sum_{\alpha = 1}^N g^{ij} \frac{\partial u^\alpha}{\partial x^i} \frac{\partial u^\alpha}{\partial x^j}.\]

There are several reasons that underpin the interest on this minimization problem. To start with,  it corresponds to the Rudin-Osher-Fatemi model \cite{ROF}, which is well-established for image denoising in the Euclidean scalar case. In fact,  given $f: \Omega \subset \R^2 \rightarrow \R$ an observed noisy grayscale image, where $\Omega$ is a bounded domain, which may be thought of as the computer screen, the goal is to reconstruct the undistorted  or ideal image $u: \Omega \rightarrow \R$ by means of the additive ansatz $f = u \, + \, \eta,$
where $\eta(x)$ are independent identically distributed zero-mean Gaussian random variables for each $x \in \Omega$.

The minimization of $\mathcal E$ is equivalent to maximize the probability $P(u|f) \propto P(u) P(\eta)$ of recovering $u$ for known $f$, when one chooses the total variation (i.e., the first integral in $\mathcal E$) to model $P(u)$. The latter should be understood as a mild regularization term, as the look for minimizers within the space of bounded variation (BV) functions allows for discontinuities across curves, which is desirable to preserve edges or boundaries of objects in images. In turn, the second integral in $\mathcal E$ plays a fidelity role as it penalizes the restored image to be too far away from the observed one, while $\lambda$ may be regarded as a weight balancing the relative relevance of the two terms involved.

In this framework, we wonder whether regular parts of the image survive to this denoising process, that is, if the regularity of the source $f$ is inherited by the minimizer $u$. This can be obtained under natural curvature restrictions on $\mathcal N$ and topological ones on the range of $f$. More precisely, let us consider the following upper bound for the convexity radius:
\begin{equation} \label{convexity-rad}
R_\kappa := \left\{\begin{array}{cc}\frac{1}{2}\min\left\{\mathrm{inj}_p \mathcal{N}, \frac{\pi}{\sqrt{\kappa}}\right\} & {\rm if \ }\kappa>0 \smallskip \\ \frac{{\mathrm{ inj}}_p\mathcal N}{2} & {\rm if \ }\kappa\leq 0.\end{array}\right.
\end{equation}
where $\kappa$ is an upper bound for all sectional curvatures of the geodesic ball $B_h(p,R)$ for $R < R_\kappa$, and ${\rm inj}_p(\mathcal N)$ denotes the injectivity radius of the ball at its center $p$.   Recall that  $B_h(p,R)$ is strongly convex provided that $R \leq R_\kappa$ (cf. \cite[Theorem IX.6.1]{Chav}).  With this notation, we prove

\begin{thm} \label{main-thm1}
	Suppose that $f(\overline\Sigma) \subset B_h(p, R)$ for some $R < R_\kappa$. Then there exists  a  minimizer $u$ of the functional $\mathcal E$, which also satisfies $u(\overline\Sigma) \subset B_h(p, R)$.  If we further assume that $\kappa \leq 0$, then
	\begin{enumerate}
	\item[{\rm (a)}] $\mathcal E$ is geodesically convex, and hence the minimizer is unique.
	\item[{\rm (b)}] If $f \in C^{0,1}(\overline\Sigma; \mathcal N)$, then $u$ inherits the same Lipschitz regularity. 	
	 \end{enumerate}
\end{thm}
\noindent	Notice that if $\mathcal N$ is simply connected, the classical Cartan-Hadamard theorem ensures that it is diffeomorphic to $\R^n$, and thus any ball of finite radius in $\mathcal N$ fulfils the small range condition for $f(\overline \Sigma)$ in the above statement. Let us further remark that the non-positivity of the curvature plays a crucial role in (a) since one can construct counterexamples where the total variation fails to be convex in case $\mathcal N = \mathbb S^2$ for 1-dimensional domains (see \cite{GLM?}). Concerning (b), we highlight that Lipschitz continuity is the optimal regularity that one could expect even if the source $f$ is smooth and with compact support.

As in the context of image processing the map $f$ has typically jumps, it is important to know which regular parts of the image remain after denoising. Therefore let us stress that a local version of our main result also holds. Indeed,
\begin{cor} \label{cor2}
	Assume that $\mathcal N$ is non-positively curved and $f(\overline \Sigma)$ is contained in a strongly convex geodesic ball of $\mathcal N$.  If we further assume that $f \in C^{0,1}(\overline \Sigma; \mathcal N)$ locally, then $u$ is also locally Lipschitz in $\overline\Sigma$.	
\end{cor}

We also prove a stronger regularity result for manifold-valued signal denoising, that is, when the domain is 1-dimensional. More precisely,
\begin{thm} \label{thm:1dcase}
	Let
$f: \Gamma \rightarrow \mathcal N$, where $\Gamma$ is either an interval $[0,1]$ or $\mathbb S^1$, so that $f(\Gamma)\subset B(p,R)$ with $R<\frac{{\rm inj}_p\mathcal N}{2}$. Then, there exists a minimizer $u$ of $\mathcal{E}$. Moreover,
\begin{itemize}
\item[(a)] If  $f \in C^{0,1}(\Gamma; \mathcal N)$, then $u\in C^{0,1}(\Gamma;\mathcal N)$.
\item[(b)] If $f\in BV(\Gamma;\mathcal N)$ and $R<\min\left\{\frac{\mathrm{inj}_p \mathcal{N}}{2}, \frac{\pi}{4 \sqrt{\kappa}}\right\}$ then \begin{equation}\label{local_estimate}
                                                              |u'|\leq C |f'|\,,\quad {\rm for\ some
  \ } C>0\,, {\rm \ as \ Borel\  measures.}
                                                            \end{equation}
\end{itemize}
\end{thm}
Let us highlight that item (b) is a fully local estimate for a solution of a variational problem, and to the best of our knowledge this is the first instance of such a result for the elliptic manifold-valued case (see \cite{GiacomelliLasicaMoll2019} for the parabolic case), which was an open problem posed in \cite{GL19}, since generic Riemannian targets are the natural setting for one-dimensional multichannel data. Notice that the one-dimensional domain gives such a stronger conclusion under milder assumptions (no curvature restriction and BV source).

%

 In addition, the proof of Theorem \ref{main-thm1} exploits a regularized version of the ROF-functional, which leads to a perturbed 1-Laplacian problem of independent interest because of its relation with fluid mechanics (see Theorem \ref{main-thm3}). Before entering into details of this further outcome, let us contextualize our main results and why they amount to several improvements to the current state of the art.

\subsection{Main results in the perspective of the literature}

Despite the plethora of applied examples in computer vision (see Section \ref{manifold-rel} for details) where the Euclidean space is not the best match to encode complicated real data, to the best of our knowledge, there is no previous study of a ROF-like model neither in the manifold-constrained setting nor on curved surface domains. In fact, not even the regularity for the vectorial version of the problem (that is, for $\mathcal N = \R^n$) is tackled, due to the lack of classical results up to the boundary.

In contrast, the  problem  has been widely studied in the scalar Euclidean case $f: \Omega \subset \R^2 \rightarrow \R$ (cf.~\cite{CCN07, CCN11}, or \cite{CCN15, CCNS} for detailed surveys). Roughly speaking, Caselles, Chambolle and Novaga showed in \cite{CCN07} that the ROF model creates no new discontinuities in addition to those already present in the observed $f$. Later on, the same authors proved in \cite{CCN11} that local Hölder regularity of the data is acquired by the minimizer, with the same exponent. Moreover, if the domain $\Omega$ is convex, they showed that the same is true for global Lipschitz regularity. The extra convexity restriction was removed for the continuous case by Mercier in \cite{Merc}, and the general Lipschitz situation is addressed by Porretta in \cite{Porr} thanks to a Bernstein-type argument, which we successfully adapt here to curved spaces both in the domain and the target.

In this framework, our main theorem deals with the corresponding Lipschitz regularity problem for vector-valued functions $u: \Omega \rightarrow \R^N$ which minimize the functional $\mathcal E$ given by \eqref{eq:rof-functional} subject to a family of smooth non-linear (manifold) constraints. In this sense, our results can be regarded as the manifold constrained counterpart of \cite{CCN11, Porr}, with the additional technical difficulty of considering also non-flat surfaces $\Sigma$ (abstract 2-manifolds) as domains. Such a mini\-mization problem leads to the following non-linear elliptic system with Neumann boundary conditions:
\begin{equation}\label{EL-ROF-intr}
	\left\{ \begin{array}{rcll}
		\ds	\tau_1(u):=\dv_g \bigg(\frac{du}{|d u|}\bigg)  & = & -\lambda\exp^{-1}_u \!f &\qquad \text{in $\Sigma$}, \medskip \\
		\nu\cdot { d u}& = & 0 &\qquad \text{on $\partial\Sigma$},
	\end{array}\right.
\end{equation}
 The precise definition of the function $\frac{du}{|du|}$ and the Neumann condition will be given in Section \ref{tensions} and Definition \ref{eq:rof-sol}.

The regularity of the system above appears to be much harder than the parabolic counterpart (see \cite{GiacomelliLasicaMoll2019} for $\lambda = 0$), as evidenced by the absence of previous references. In short, we generalize the existing literature in three different directions: vector-valued, manifold-constrained and curved domain.

\subsection{A geometric viewpoint to ground the naturality of the hypotheses} Apart from the already mentioned interest coming from image processing, the minimization of \eqref{eq:rof-functional} is of geometric interest because the corresponding Euler-Lagrange equations \eqref{EL-ROF-intr} can be interpreted in the scalar case as the following prescribed mean curvature problem:
\[H = -\lambda \exp^{-1}_t \!f=\lambda(t-f) \qquad \text{on} \quad \partial E_t,\]
where $H$ represents the mean curvature of the boundary of the sup-level sets $E_t:=\{u > t\}$. In fact, the sought regularity statement can be interpreted in terms of distance between any two level sets $E_t$ and $E_s$ of the minimizer $u$. This connection is unsurprising if we recall that the total variation, when applied to a characteristic function $u = \chi_E$ is the perimeter of $E$ in $\Sigma$.

In a somewhat different spirit, \eqref{EL-ROF-intr} can be regarded as a shifted spectral problem for the so-called 1-tension or rough 1-Laplacian (see Section \ref{tensions} for technicalities of the definition) associated with the eigenvalue $\lambda$ and with respect to a model map $f$, which was introduced by Jost and Kourouma in \cite{Jost-K} for the standard 2-tension or rough Laplacian.

In this sense, the $\lambda = 0$ case corresponds to 1-harmonic maps, and hence our results can be regarded as the  extension to $p =1$ of regularity theory of $p$-harmonic maps, which is extensively studied mainly for $1 < p < \infty$, and which is itself an extension of the classical harmonic map theory ($p = 2$)  started by Eells and Sampson \cite{EllS} or Schoen and Uhlenbeck \cite[Corollary, p. 310]{SchUh}, where the non-positive curvature (NPC) hypothesis and the small range condition pop up naturally. Indeed, it is well-known (see e.g. examples in \cite{Riv}) that there is no hope of getting full regularity results unless we impose topologic and/or geometric restrictions on the target. Notice that our $p =1$ case is analytically much harder because the anisotropic diffusivity of the 1-tension, meaning that its ellipticity degenerates in the direction of $d u$, while it becomes singular in directions where $u$ is constant.

Indeed, our small range condition was first considered for $p = 2$ by Hildebrandt, Kaul and Widman \cite{HKW} to prove regularity results, which were extended to $p\geq 2$ and $1< p<2$ in \cite{FaBo} and \cite{GuoXi}, respectively. As all these papers deal with the friendlier scenario of Dirichlet boundary conditions, our results are closer to the Neumann problem for harmonic maps studied by Hamilton in the NPC case \cite{Ham}. Let us point out that the lack of any harmonicity hypothesis is a constant technical difficulty in our work, because we have to control tension terms overall which in these previous references did not play any role because they vanish.

\subsection{About the meaningfulness of the generalization to manifolds } \label{manifold-rel}

Let us highlight that the interest of the manifold constrained version of the ROF model is not by the sake of meaningless generalization/abstraction; indeed, there are several instances in applied sciences for which image processing requires this setup (see e.g. \cite{SIAM-2, SIAM-1}). The easiest example to grasp the need of manifold targets is the tracking of images needed for animation or video surveillance, where data take values in the space of rigid motions of $\R^3$.

Having applications in mind, relevant examples in biomedical imaging are constrained to NPC manifolds, like DT-MRI  (Diffusion tensor magnetic resonance imaging) used in neurology, where six different data sets measured with different magnetic fields to govern the diffusion of water molecules is considered for each pixel, and this information is encoded in the space of
symmetric positive definite $3\times 3$ matrices. In a completely different direction, hyperbolic space has become specially relevant in image segmentation \cite{HypIS} and machine learning \cite{survey}, where training of the model adds some noise to the process and denoising algorithms may lead to better pattern recognition. An intuitive reason behind this is that in hierarchies (such as trees or taxonomies)  the number of nodes increases exponentially, which mimics the volume growth of a ball in hyperbolic ambients,  making the latter better suited than Euclidean traditional geometry, where the growth is polynomial.

On the other hand, positively curved targets (specially round spheres) also play a role in RGB image processing \cite{TSC} or in  EBSD (electron backscatter diffraction) methods used in materials science to study the microscopic structure of polycrystals, where data take values in ${\rm SO}(3)$ (see \cite{micro, M3D}). In turn, curved domains arise naturally in sight of non-flat computer screens, geophysical and medical applications (wind directions or cortical surfaces in brain images) or 3D colouring of surfaces for art or design \cite{Lai}.

In short, dealing with manifold-constrained signals and images is a new challenge coming from real problems that has been already addressed in the applied literature (see e.g. \cite{recent}), which can be considered in this setting some steps ahead the analytical foundations, which are somehow missing. This motivates the present paper.

\subsection{Steady Mosolov problem }

To prove existence of minimizers by means of the direct method in the Calculus of Variations, we modify the functional in \eqref{eq:rof-functional-0} by sticking a parameter $\sigma$ and introducing a Dirichlet integral to gain some coercivity, that is, we define a new energy as
\begin{equation}\label{eq:rof-functional}
\mathcal E_\sigma(u) := \int_\Sigma |d u| \,d\mu_g + \frac{\lambda}{2}\int_\Sigma \dg^2(u, f)d\mu_g + \frac{\sigma}{2}\int_\Sigma |d u|^2 \,d\mu_g,
\end{equation}
so that $\mathcal E_0 = \mathcal E$.  The corresponding Euler-Lagrange equations lead to
\begin{equation} \label{Mosolov}
\tau_1(u) + \sigma \, \tau(u) = -\lambda\exp^{-1}_u\! f \quad \text{in} \quad \Sigma,
\end{equation}
coupled with Neumann boundary conditions, where $\tau(u)$ stands for the usual tension or rough Laplacian. This is known as perturbed 1-Laplacian system, where the perturbation allows to expect for better regularity properties (see some Euclidean results in this direction, but without boundary conditions, in \cite{GiTsu} for convex scalar solutions, and \cite{Tsu} in the vector-valued situation, removing convexity). This intuitive expectation comes from the extra diffusivity (and thus ellipticity) coming from the usual tension, and in our case allows for regularity results even for $\kappa > 0$. Indeed, we prove an enhanced version of Theorem \ref{main-thm1} for minimizers of $\mathcal E_\sigma$.

\begin{thm} \label{main-thm3}
		Suppose that $f(\overline\Sigma) \subset B_h(p, R)$ for some $R < R_\kappa$. Then for any fixed $\sigma > 0$, there exists a  minimizer of the functional $\mathcal E_\sigma$, which also satisfies $u(\overline\Sigma) \subset B_h(p, R)$.
		\begin{enumerate}
			\item[{\rm (a)}] If we further assume that $\kappa \leq 0$, then $\mathcal E_\sigma$ is geodesically convex, and hence the minimizer is unique. If $f \in C^{0,1}(\overline\Sigma; \mathcal N)$, then $u$ inherits the same Lipschitz regularity.
			\item[{\rm (b)}] In the case $\kappa > 0$ there exists a minimizer of $\mathcal E_\sigma$ which is also Lipschitz continuous, provided that $f(\overline\Sigma) \subset B_h(p, R)$ for some $R < \min\left\{\frac{\mathrm{inj}_p \mathcal{N}}{2}, \frac{\pi}{4 \sqrt{\kappa}}\right\}$.
		\end{enumerate}
	\end{thm}
Notice that the stronger small range condition in (b) was also required in \cite{Fuchs} to prove regularity of stationary $p$-harmonic maps for $p \geq 2$.

Let us highlight that the sought of minimizers for this functional is of independent interest since, if we view the vector in the right hand side of \eqref{Mosolov} as the gradient of some pressure function, then this equation models an isothermal
steady laminar flow of a  non-Newtonian (Bingham) incompressible fluid, where $u$ stands for the flow velocity and $\Sigma$ is the cross-section of and infinitely long duct. In this context, our Neumann condition means that there is no flux across the boundary, and the parameter $\sigma$ controls the plastic viscosity. This is usually referred to as Mosolov problem \cite{Mos}.

\subsection{Strategy for the proof of Lipschitz regularity} Let us give a brief sketch of the proof for Lipschitz regularity. To remove the degeneracy of our system, we need to stick in \eqref{eq:rof-functional} a further parameter $\varepsilon > 0$ to regularize the energy density in the first term by shifting it to $\sqrt{|d u|^2 + \varepsilon^2}$. Now the idea is to apply a Bernstein type technique, that is, to differentiate the resulting PDE system
\[	\dv_g \bigg[\Big(\frac1{\sqrt{|du|^2 + \varepsilon^2}} + \sigma \Big) du\bigg]  =  -\lambda\exp^{-1}_u f \qquad \text{in $\Sigma$}\]
in order to prove that $w = |du|^2$ is a subsolution of an elliptic problem to which one can apply maximum principle arguments to get pointwise estimates at a maximum point. Just to begin with this strategy we have to deal with two technical difficulties not present in the previous literature: the right hand side is only Lipschitz and, unlike the scalar Euclidean case from \cite{Porr}, solutions of this regularized problem are not granted to be smooth by classical standard theory.

To overcome this issue we substitute $f$ with a smoothed version of it constructed via distance mollification as in \cite{Karcher1997} (cf. Section \ref{boost}) and we need to show that solutions of the regularized problem are at least $C^3$ up to the boundary, and this is precisely a delicate point, as regularity theory for Neumann elliptic problems is not at hand, and we have to adapt techniques from \cite{Arkhipova} to our particular setting in order to prove Hölder estimates for the gradient of $u$ almost from scratch. Moreover, even interior regularity is only guaranteed in case $\sigma > 2 C(R) R$ (cf. \cite{GiMo, HW}), but we need to pass to the limit $\sigma \to 0$ afterwards.

 A crucial issue is that the arguments only work for systems in the Euclidean space, and hence we need to extend our problem suitably to a tubular neighborhood of the target manifold, following ideas from our previous work in \cite{CRMP}, see details in Section \ref{setupA}. After the extension and checking that  we get a system of natural growth on the normal part, we can readjust our results in  \cite{CRMP} to get an intrinsic Caccioppoli inequality and higher integrability estimates for solutions of the present problem (see Proposition \ref{prop:CMP-inequalities} for details).

To achieve Hölder regularity, we start by flattening the boundary covering it by a system of isothermal coordinates. Then the idea is to freeze a point $\bar x$ in the domain and introduce an auxiliary function $v$ which solves an homogeneous system with mixed boundary conditions in a suitable neighbourhood of $\bar x$ and for which regularity is known. Next, after a chain of estimates for the difference $v - u$ performed by choosing appropriate test functions, we prove that $u \in C^{0, \beta}(\overline \Sigma; \mathcal N)$ for any $\beta \in (0, 1)$ (see Theorem \ref{thm:holder-u}). Working harder to get finer estimates, we reach in Theorem 	\ref{thm:holder-Du} that $u\in C^{1,\beta_0}(\overline{\Sigma}; \mathcal N)$ for some $\beta_0 \in (0,1)$. The solutions are shown to be actually smooth up to the boundary (see Proposition \ref{prop:holder-D2u}) by combining Schauder estimates and bootstrapping.

It is a crucial point here that the above arguments can only succeed by working in extrinsic coordinates, where there is PDE theory available that can be accommodated after non trivial adaptations. Additionally, the consideration of a curved surface instead of a flat domain imposes extra complications that cannot be neglected by covering $\overline \Sigma$ with coordinate neighbourhoods and arguing as in a flat domain. Indeed, precisely the non-Euclidean metric yields $x$-dependence within the elliptic operator that governs the extended problem, which certainly makes the proof of the Hölder estimates  more intricate than the corresponding version for a subset of $\mathbb R^2$.

Notice that the aforementioned smoothness is not enough for our purposes, because the bounds arising from the proofs are not uniform either in $\sigma$ or $\varepsilon$, and hence the estimates will not survive after passing to the limit to get regularity of the original ROF problem. But this process allows to differentiate the equation and proceed with the Bernstein technique (see Section \ref{Lip}).
 Surprisingly enough and in contrast with the method we have just described, the latter can only work by working intrinsically (without using Nash embedding). This is why the interplay between analytical and geometric approaches is essential to achieve our main results.

\subsection{Structure of the paper} The paper is structured as follows. In Section \ref{notation} we recall the extrinsic versus intrinsic way to look at the $p$-Laplacian problems and Sobolev spaces on manifold, we also fix notation, conventions and definitions that will be repeatedly used throughout the whole paper. Section \ref{sec:tv_def} includes the technicalities that are hidden on purpose for the sake of simplicity of the introduction to define the minimization problem on the right space of BV functions between manifolds, which also amounts to define an appropriate notion of solution and prove lower semicontinuity of the corresponding relaxed energy functionals (cf. Theorem \ref{tv-relax}), which will be crucial for the existence of minimizers.

In turn, uniqueness follows by direct methods of calculus of variations, provided the functional is convex. Next, Section \ref{convex} is devoted to prove geodesic convexity for an even more general energy functional, where we have to derive the second variation formula from scratch, since in all the previous literature both tension and boundary terms are neglected by harmonicity and boundary conditions, respectively. This is used to derive convexity for smooth functions between Riemannian manifolds in the NPC case. To conclude the argument, one needs a density of the latter within Sobolev functions, which is not true without serious topological restrictions; to overcome this, we apply an idea from \cite{PiVe} of flattening out the ends outside a geodesic ball to construct a global system of normal coordinates where the squared distance function is still convex (see Lemma \ref{lemma:density}). As a by-product, we achieve geodesic convexity for the total variation (cf. Theorem \ref{convex-tv}). Then we are in position to prove existence of minimizers in Section \ref{exist}, after showing that the composition of a map $u$ with a retraction into the geodesic ball decreases the energy $\mathcal E_\sigma$ and hence the small range condition also holds for minimizers (see Proposition \ref{prop:comparison}).

 The aforementioned arguments for Hölder and Lipschitz regularity are carried over in Sections \ref{Hol-reg} and \ref{Lip}, respectively. Let us remark that the Lipschitz regularity with a uniform bound for the regularized problem depending on $\varepsilon, \delta$ follows directly from the elliptic estimate for $w =|d u|^2$ (see Proposition \ref{Porr-cvx}) in the convex case. In order to remove the latter extra assumption, we need to multiply $w$ by a suitable test function, which only depends on the domain surface, so that one can still apply Hopf's maximum principle (see Section \ref{peso}). In this way, one can prove Lipschitz regularity for the approximated problem in Proposition \ref{prop:lipschitz}, and for the original one after a limiting process (cf. Proposition \ref{final_result}). The corresponding statement for the Mosolov problem is more complicated, as it requires multiplication of $w$ by a test function that depends both of $\Sigma$ and the target manifold, cleverly chosen so that the quadratic terms in $w$ have a good sign which allows to conclude via maximum principle (see Section \ref{LipMos}).

Concerning the case of signal denoising, in Section \ref{signal} we indicate how to get the Lipschitz regularity via simpler arguments than in the previous sections. The local estimate for BV functions is proved by means of local integral estimates and an approximating scheme similar to those carried out in \cite{GL19} for the vectorial case, but properly enhanced so that the arguments apply to manifold constrained settings.
	
	\section{Notation, conventions and background material} \label{notation}

	Let us be a bit more precise about the definition of relevant operators and Sobolev spaces for maps between manifolds.
	
	\subsection{$p$-tension or rough $p$-Laplacian (intrinsic versus extrinsic viewpoint)} \label{tensions}
	
	Let $(\mathcal N, h)$ be as before, and take $(\mathcal M, g)$ be an $m$-dimensional compact Riemannian manifold with smooth boundary $\partial \mathcal M \neq \emptyset$. Recall that for each point $x \in \mathcal M$ the differential of $u$ at $x$ is a linear map from $T_x \mathcal M$ to $T_{u(x)}\cn$, namely, it holds $d u(x) \in {\rm Hom}(T_x \mathcal M, T_{u(x)}\cn)$, and the latter space is in turn isomorphic to the tensor product $T_x^\ast \mathcal M \otimes T_{u(x)}\cn$.
	
	We take $(x^1, \ldots, x^m)$ and $(y^1, \ldots, y^n)$ local coordinates in $\mathcal M$ and $\mathcal N$, respectively. Then $du$ as a section of the bundle $\mathscr B:=T^\ast \mathcal M \otimes u^\ast T \mathcal N$ can be expressed locally as
	\[d u = \partial_i u^\alpha \, dx^i \otimes \frac{\partial}{\partial y^\alpha}\bigg|_u, \qquad \text{with} \quad \partial_i:=\frac{\partial }{\partial x^i }.\]
	For convenience and brevity, we will combine $du$ with the notation $u_\ast$ as both are quite usual in the literature. Hereafter we use the convention that an index repeated as sub and superscript in
	a product means summation over the range of the index, where Latin indices range in $\{1,\ldots, m\}$, while Greek letters take values on $\{1,\ldots, n\}$.
	
	From an intrinsic viewpoint, given $p\geq 1$, the $p$-Dirichlet energy (or total variation functional for $p = 1$) is given by
	\begin{equation}\label{p-energy}\mathscr E_p(u) = \frac1{p} \int_{\mathcal M} |d u|^p \, d\mu_g,\end{equation}
	where $d\mu_g$ denotes the Riemannian volume element and the energy density can be written as
	\[|d u|^2 = \<d u, du\>_{g \boxtimes h} := \tr_g(u^\ast h) = g^{ij} (h_{\alpha \beta} \circ u) \partial_i u^\alpha \partial_j u^\beta, \]
which follows from the definition of the pull-back metric $(u^\ast h)_{ij} :=h(\partial_i u, \partial_j u)$.

 Given  $\nu$ the unit outward normal to $\partial\M$, we interpret $\nu\cdot du$ as the inner product with respect to the induced metric on the boundary. Moreover, notice that $\nu \cdot du$ is actually an abuse of notation for $du(\nu) =u_\ast \nu \in T_u \mathcal N$; for convenience, sometimes we will also use the notations $\nabla_\nu u$ or $\frac{\partial u}{\partial \nu}$.

On the other hand, let $\nabla$ be the Levi-Civita connection  of $(\mathcal M,g)$ and $^{h}\nabla$ the one of $(\mathcal N,h)$. Set $\widetilde\nabla$ the induced connection on $u^* T\mathcal N$, which is defined as follows: for $X\in T\mathcal M$, $V\in u^* T\mathcal N$,
\begin{equation} \label{ind-conn}
\widetilde\nabla_X V:=^{h}\!\!\nabla_{u_* X} V\in C^\infty(u^*T\mathcal N).
\end{equation}
In turn, if $\nabla^\ast$ denotes the connection on $T^\ast \mathcal M$ dual to $\nabla$, we will work with $\nabla^{g \boxtimes h}:= \nabla^{\ast} \otimes \widetilde{\nabla}$, which is the unique linear connection of $\mathscr B$ so that
\[\nabla^{g \boxtimes h}(A \otimes B) = (\nabla^\ast A) \otimes B + A \otimes \widetilde{\nabla} B \]
for $A \in C^\infty(T^\ast \mathcal M)$ and $B \in C^\infty(u^\ast T\mathcal N)$. Additionally, such a connection fulfils the compatibility condition
\[v\<X, Y\>_{g\boxtimes h} = \big<\nabla^{g \boxtimes h}_v X, Y\big>_{g\boxtimes h} + \big<X, \nabla^{g \boxtimes h}_v Y\big>_{g\boxtimes h}\]
for all sections $X, Y \in C^\infty(\mathscr B)$ and $v \in T\mathcal M$.

	The critical points of $\mathscr E_p$ are known as $p$-harmonic maps, and have vanishing $p$-tension or rough $p$-Laplacian, that is, satisfy $\tau_p(u) =0$, where $\tau_p(u):=\tr_g\big(\nabla^{g \boxtimes h} du\big) = \tau_p(u)^\gamma \frac{\partial}{\partial y^\gamma}\Big|_u \in C^\infty(u^\ast T \mathcal N)$, which in local coordinates reads as
\begin{align*}
\tau_p(u)^\ga & =\dive_g\big(|du|^{p-2} d(u^\gamma)\big) + |du|^{p-2} g^{ij} \, ^{h}\Gamma_{\alpha \beta}^\gamma(u) \partial_i u^\alpha \partial_j u^\beta
\\ & = \frac1{\sqrt{\det g}} \partial_i\Big(|du|^{p-2} \sqrt{\det g} \, g^{ij} \partial_j u^\gamma\Big) + |du|^{p-2} g^{ij} \, ^{h}\Gamma_{\alpha \beta}^\gamma(u) \partial_i u^\alpha \partial_j u^\beta
\end{align*}
for $\gamma =1, \ldots, n$, where $^{h}\Gamma_{\alpha \beta}^\gamma(u)$ denote the Christoffel symbols for the Levi-Civita connection of the metric $h$. In particular, we will deal with the cases $p =1,2$, that is,
\begin{align*}
\tau(u)^\gamma = \tau_2(u)^\gamma = \Delta_g(u^\gamma) + g^{ij} \, ^{h}\Gamma_{\alpha \beta}^\gamma(u) \partial_i u^\alpha \partial_j u^\beta =  g^{ij}\Big(\partial^2_{ij} u^\gamma - \, ^{g}\Gamma_{ij}^k \partial_k u^\gamma + \, ^{h}\Gamma_{\alpha \beta}^\gamma(u) \partial_i u^\alpha \partial_j u^\beta\Big)
\end{align*}	
and, with a shorter notation,
\begin{equation} \label{tau1-int}
\tau_1(u) = \dive_g\bigg(\frac{d u}{|d u|}\bigg) + \frac1{|d u|} \Gamma_u (du, du).
\end{equation}

Alternatively, from an extrinsic perspective, Nash theorem \cite{Nash1956} guarantees that there exists an isometric embedding $\iota: \mathcal N  \hookrightarrow \R^N$ and, by completeness of $\mathcal N$, we even know by \cite{OM} that $\iota(\mathcal N)$ is a closed subset of $\R^N$ (in other words, that the embedding is proper). Then we can work with the extended map
\[\breve{u}:= \iota \circ u: \mathcal M \longrightarrow \R^N.\]
In this setting, with $\nabla \breve u =(\nabla \breve u^1, \ldots, \nabla \breve u^N) \in \R^{m N}$, it holds
\[|\nabla \breve u|^2 = \sum_{\alpha = 1}^N g^{ij} \partial_i \breve u^\alpha \, \partial_j \breve u^\alpha \quad \text{and} \quad \tau_p(\breve{u}) = \Big(\dive_g \big(|\nabla \breve u|^{p-2} \, \nabla \breve u\big)\Big)^\top,\]
where $\R^N \ni X \mapsto X^\top \in T_u\mathcal N$ denotes the orthogonal projection onto the tangent space. As the normal component can be represented by means of the second fundamental form $\mathcal A_u$ of $\mathcal N$, we get the following expression:
\[\tau_p(\breve{u})^\gamma = \dive_g \big(|\nabla \breve u|^{p-2} \, \nabla \breve u^\gamma\big)\Big) + |\nabla \breve u|^{p-2} g^{ij} \mathcal A_u^\gamma(\partial_i \breve u, \partial_j \breve u), \qquad 1 \leq \gamma \leq N.\]
being $\mathcal A_u(X, Y) \in (T_u \mathcal N)^\perp$ the second fundamental form of $\mathcal N$ given by $\mathcal A(X, Y) = -(D_X Y)^\perp$, where $D$ is the standard directional derivative in $\R^N$ and $X, Y$ are extended arbitrarily as tangent vector fields to $\mathcal N$ in a
neighborhood of $u\in \mathcal N$. Notice that, with a shorter notation for $p =1$, we have
\[\tau_1(\breve u) = \dive_g\bigg(\frac{\nabla \breve u}{|\nabla \breve u|}\bigg) + \frac1{|\nabla \breve u|} \mathcal A_u (\nabla \breve u, \nabla \breve u),\]
which formally resembles its intrinsic counterpart \eqref{tau1-int}. Hereafter, for simplicity, we will always write $u$ instead of $\breve u$, unless the meaning is unclear from the context.

\subsection{$p$-energy on surfaces}	\label{p-energy-sect}
In the case that $\mathcal M=\Sigma$, we can further use isothermal coordinates on the surface. Indeed, let $\big\{(U_\ell,\phi_\ell)\big\}$ be an atlas of $\Sigma$ such that the metric $g$ within $U_\ell$ can be written as $g_{ij}=\varrho_\ell^2 \, \delta_{ij}$ with $\varrho_\ell \in C^1(U_\ell;\R^+)$. 
In these coordinates, set $\Omega_\ell:=\phi_\ell(U_\ell)$, and thus the $p$-energy functional \eqref{p-energy} can be written as
\begin{equation}\label{p-energy-extrinsic}
\mathscr E_p(u) = \frac{1}{p}\sum_{\ell}\int_{\Omega_\ell} \big(\chi_\ell \cdot \varrho_\ell^{2-p} |\nabla u|^p\big)\circ\phi_\ell^{-1}(x) \ dx =: \frac{1}{p}\int_{\Omega}\varrho^{2-p}|\nabla u|^p\,dx,
\end{equation}
with $\{\chi_\ell\}_\ell$ a partition of unity subordinate to $U_\ell$. Notice that, for simplicity, we just write $\varrho_\ell$ instead of $\chi_\ell \cdot \varrho_\ell$ and, as usual, we omit composition with $\phi_\ell^{-1}$, unless the meaning is unclear from the context.


Now, observe that the existence of isothermal coordinates up to the boundary (e.g.~\cite[Lemma 4]{YZ}) permits to take a finite covering $\{U_\ell\}_{\ell\in I}$ and to suppose that there exists a constant $C_p>0$ such that
\begin{equation}\label{boundforiso}C_p\leq \rho_i^{2-p}\big|_{\Omega_\ell} \leq \frac{1}{C_p}\, \quad {\rm for \ all \ }\ell\in I {\rm \  and \ } p\geq 1.\end{equation}

\subsection{Extrinsic versus intrinsic Sobolev spaces}
Given a smooth compact manifold $\mathcal M$, with possibly smooth boundary $\partial\M$,  $L^p(\mathcal M)=L^p(\mathcal M;\R^N)$  ($1\leq p < \infty$) and $L^\infty(\mathcal M)=L^\infty(\mathcal M;\R^N)$ denote, respectively,  the set of Lebesgue integrable functions such that the $p$-power of the modulus is Lebesgue integrable in $\mathcal M$ and the set of essentially bounded measurable functions in $\mathcal M$. Complementarily, for $p\geq 1$ the spaces $W^{k,p}(\mathcal M)=W^{k,p}(\mathcal M;\R^N)$ denote  the sets of functions in $L^p(\mathcal M)$  such that its derivatives up to order $k$ have finite $L^p$ norm. In particular, we denote by $H^1(\mathcal M)$ the space $W^{1,2}(\mathcal M)$. We recall that $W^{1,p}(\mathcal M)$ is the closure of $C^\infty_c(\mathcal M)$ functions in the standard $W^{1,p}$-norm.

	Let $\Omega\subset\R^m$ be an open bounded set with Lipschitz boundary $\partial\Omega$ and  $\Gamma \subset \partial\Omega$. For \(1 \leq p < m\), \(W^{1,p}_\Gamma(\Omega)\) is the subset of \(W^{1,p}(\Omega)\) consisting of functions \(u \in W^{1,p}(\Omega)\) with \(Tu = 0\) on \(\Gamma\), where \(T\) is the trace operator on \(\partial \Omega\). In particular, when \(\Gamma = \partial \Omega\), \(W^{1,p}_\Gamma(\Omega)\) is denoted as \(W^{1,p}_0(\Omega)\).

	We will use the following generalization of the Poincar\'e inequality. We omit the proof since it is verbatim that for the standard $W^{1,p}_0(\Omega)$ case (e.g., \cite[Thm 3, p. 279]{Evans}).
	
	\begin{lemma}\label{lem:poincare}
		Let \(u \in W^{1,p}_\Gamma(\Omega)\) for some \(1 \leq p < m\). Then, the inequality
		\begin{equation*}
		\lVert u \rVert_{L^q(\Omega)} \leq C \lVert Du \rVert_{L^p(\Omega)}
		\end{equation*}
		holds for each \(q \in [1, \frac{mp}{m-p}]\), where \(C > 0\) is a constant depending only on $p,q,m$ and $\Omega$.
	\end{lemma}

Concerning $\mathcal N$-valued functions, there are two different notions of Sobolev spaces: an extrinsic version and an intrinsic one. For the first, identifying again a function $u$ with its extension $\iota\circ u$,
\begin{equation} \label{Sob-ext}
W^{1,p}(\mathcal M;\mathcal N):=\{u\in W^{1,p}(\mathcal M) : u(x)\in \mathcal N  {\rm \ for \  a.e. \   }x \in \mathcal M \}.
\end{equation}
It turns out (\cite[Proposition 2.7]{convent-van}) that this definition does not depend on the choice of the embedding $\iota$, contrary to the case of higher order Sobolev spaces (see \cite{convent-van-2}).

On the other hand, in the event that the $n$-dimensional manifold $\mathcal N$ is covered by a normal coordinate chart centered at some point $p\in\mathcal N$, the manifold $\mathcal N$ can be identified with the Euclidean space $\R^n$ endowed with a metric that can be expressed globally in normal Cartesian coordinates. Identifying thus $u:\mathcal M\to\mathcal N$ with its vector valued representation, one can define the intrinsic Sobolev space of first order as
\begin{equation} \label{Sob-in}
W_{\rm in}^{1,p}(\mathcal M;\mathcal N):=W^{1,p}(\mathcal M).
\end{equation}

\subsection{Further notation and conventions}\label{sec:notations}

%
%
We now recall a mollification technique given in \cite{Karcher1997} for Lipschitz continuous functions, $C^{0,1}(\mathcal M;\mathcal N)$:
 Firstly, for a given $0<\delta$ 
 smaller than the injectivity radius of $\mathcal M$, we 
 consider the mollifier $\varphi_\delta:\mathcal M\times\mathcal M \to \mathbb R$ defined as
	\[
	\varphi_\delta(x,y) = \psi\left(\frac{d_g(x,y)}{\delta}\right)\left[\int_{B(0,\delta)}\psi\left(\frac{d_g(x,y)}{\delta}\right)d \mu_{{\rm exp}^*_x g}\right]^{-1},
	\]
	where $B(0,\delta) = \exp_x (B_g(x,\delta))$ is regarded as a ball in $\R^m$, while $\psi\in C^\infty(\mathbb R)$ is such that $\psi((-\infty, 1/4]) = 1, \psi([3/4, \infty)) = 0$ and $\psi' \leq 0$. 
	Then, if we further assume that $f(\mathcal M)\subset B_h(p,R)$ for some $p\in\mathcal N$ and $R<R_\kappa$,  we can define $f_\delta: \mathcal M \to \mathcal N$ as
\begin{equation} \label{moli-f}
	f_\delta(x) := {\rm argmin }_{p\in\mathcal N} \left\{\frac{1}{2}\int_{\mathcal M} d_g^2(p,f(\cdot))\varphi_\delta(x,\cdot) d\mu_g \right\}.
	\end{equation}
 Let us point out that \cite[Theorem 2.1]{afsari} ensures that $f_\delta(x)$ is well defined and, 
%
by \cite[Theorem 4.3]{Karcher1997}, we have that $\{f_\delta\}_{\delta}\subset  C^\infty(\mathcal M; \mathcal N)$ for $\delta$ small enough. In addition, as $f\in C^{0,1}(\mathcal M;\mathcal N)$ with Lipschitz constant $L_f$, \cite[Theorem 4.4, 4.6]{Karcher1997} implies that
	$f_\delta(p) \to f(p)$ as $\delta \to 0$ and $f_\delta$ are Lipschitz functions with Lipschitz constant $C L_f$.
	
\medskip	
Given an open set $\Omega\subset\R^m$, we represent by $\mathcal L^{p,\lambda}(\Omega)$  the Campanato space, which consist of those functions $u\in L^p(\Omega)$  satisfying
	\[
	\sup_{\substack{ x\in\Omega,\\ \rho \leq \text{diam}(\Omega)}} \bigg(\rho^{-\lambda}\int_{\Omega\cap B_\rho( x)} |u - u_{x, \rho}|^p dx\bigg) < \infty,
	\]
	such that $u_{x, \rho}$ is defined as
	\[
	u_{x,\rho} := \fint_{\Omega\cap B_\rho( x)} u \ dx.
	\]

\medskip
 We specify our conventions for  constants. Indexed letters
		$C$, i.e. $C_0$, $C_1$, etc. will retain a specific meaning throughout the whole paper, while the letter $C$ denotes a generic constant, which is always allowed to change from line to line and depends on the quantities listed in the corresponding statement. In addition, if the constants depend on some specific variable (e.g. $\varepsilon$), this will be denoted by an index with the corresponding variable (e.g. $C_\varepsilon$)

	\section{Total variation functional and technical setup for the Euler-Lagrange equations}\label{sec:tv_def}
	\subsection{Manifold constrained functionals with linear growth}
In this subsection we extend the energy functional $\mathscr E_1$ for functions not necessarily in $W^{1,1}(\Sigma;\mathcal N)$. In particular, with the notation in \eqref{p-energy-extrinsic}, we define for $u\in L^1(\Sigma;\mathcal N)$,
\begin{equation}\label{TV-def}\mathcal{TV}(u):=\left\{\begin{array}
  {ll} \!\!\! \displaystyle \sum_\ell\!\left(\int_{\Omega_\ell} \!\!\!\varrho_\ell |\nabla u|\, dx+\int_{\Omega_\ell} \!\!\varrho_\ell \, d|D^c  u| \displaystyle+\int_{\Omega_\ell \cap J_{ u}} \!\!\!\varrho_\ell  d_h(u^+, u^-)\, d\mathcal H^1\!\right)\!, & {\rm  \!\!\! if \ }  u \circ \phi_\ell^{-1}\!\in BV(\Omega_\ell;\mathcal N)
  \smallskip \\ +\infty, & {\rm  \!\!\! otherwise.}
\end{array}\right.\end{equation}
 Note that, again, we are omitting composition with $\phi^{-1}_\ell$ inside the integrals. Here we are using an extrinsic concept of BV functions; i.e. for an open set $\Omega\subset\R^m$,
$$BV(\Omega;\mathcal N):=\{u\in BV(\Omega) : u(x)\in \mathcal N \ {\rm for \ a.e. \ } x\in \Omega\},$$ and $J_{u}$ denotes the jump set of $u\in BV(\Omega)$ while $D^cu$ is the Cantor part of the Radon measure $Du$. Since we do not use any properties of BV functions in our paper, we do not enter into details of the definition or properties of the space $BV(\Omega)$, but we instead refer to \cite{AmbrosioFuscoPallara2000}.
Moreover, we say that $u\in BV(\Sigma;\mathcal N)$ if  $u\circ \phi^{-1}_\ell\in BV(\Omega_\ell;\mathcal N)$  for all $\ell\in I$.

We will prove the next result, concerning functions valued on a geodesic ball $B_h(p,R)$ around a point $p\in \mathcal N$ with radius $R$ smaller than $R_\kappa$ defined as in \eqref{convexity-rad}.

\begin{thm}\label{tv-relax} Let $R< R_\kappa$. The following holds:
\begin{itemize}
\item[(a)]
  Let $\{u_k\}_k\subset W^{1,1}(\Sigma;B_h(p,R))$ for some $p\in\mathcal N$  be such that $u_k\to u$ in $L^1(\Sigma;B_h(p,R))$, then $$\mathcal{TV}(u)\leq \liminf_{k\to\infty} \mathscr E_1(u_k).$$
  \item[(b)] Given $u\in BV(\Sigma;B_h(p,R))$ for some $p\in\mathcal N$, there exists $u_k\in C^1(\Sigma;\overline{B_h(p,R)})$, such that $\mathscr E_1(u_k)\to \mathcal{TV}(u)$ as $k\to\infty$.
\end{itemize}
\end{thm}
\begin{proof}
First, we observe that it is enough to prove the result for $\Omega:=\Omega_\ell$. For simplicity of notation, we also omit composition with $\phi^{-1}_\ell$.
Under this notation, we observe that the result is very close (except for the fact of the different open and closed geodesic balls) to prove that $\mathcal TV$ is the relaxed functional with respect to the $L^1$-convergence of $$u\mapsto\left\{\begin{array}{ll}\mathscr E_1(u) & {\rm if \ } u\in C^1(\Omega; \overline{B_h(p,R)}) \smallskip \\ +\infty & {\rm otherwise}\end{array}\right..$$

In the case that the target space $\mathcal Y$ is a compact, orientable manifold without boundary with trivial homotopy group, we could apply directly the results in \cite{GuiMu} to obtain that this is the case. However, in our case $\mathcal Y=\overline{B_h(p,R)}$ has a nonempty boundary. To fix this issue, we argue as follows.

\smallskip
We take $R<R'<R_\kappa$. Since $\partial B_h(p;R')$ is geodesically convex, we can consider the double of $\mathcal R:=B_h(p;R')$, $\mathcal R^D$, which is typically constructed  by gluing together two copies of $\mathcal R$ along the identity map of the boundary, but this will not give a regular enough closed submanifold in general.  To overcome this, we can argue as in \cite[Proposition 4.1]{CRW}, that is, roughly speaking, by perturbing the metric in a small inner neighborhood of the boundary $\partial \mathcal R$ to form a cylindrical end so that the gluing is well defined and the submanifold $( \mathcal R^D, \tilde g)$ preserves its degree of regularity and the modified metric $\tilde g$ satisfies $\tilde g|_{B_g(p,R)} = g$.

Then, \cite[Corollary 6.4]{GuiMu} applies with $\mathcal Y= \mathcal R^D$ and the lower semicontinuity in $(a)$ is a direct consequence. For $(b)$, we observe that \cite[Corollary 6.4]{GuiMu} yields the existence of a sequence $\{u_k\}_k\subset C^1(\Omega; \mathcal R^D)$ such that $u_k\to u\in L^1(\Omega;\mathcal R^D)$ and $\mathscr E_1(u_k)\to \mathcal TV(u)$ as $k\to\infty$. By the strong convergence, we can take a subsequence, not relabeled, such that $u_k\to u$ a.e. Therefore, we can suppose that $u_k\in C^1(\Omega;\overline{B_h(p,R)})$.
\end{proof}

\begin{bem}\label{rem:1dcase} In the case of a one-dimensional domain $\Gamma$, Theorem \ref{tv-relax} holds without requiring any restriction on the range, see \cite[Lemmas 2.3 and 2.4]{GiacomelliLasicaMoll2019}.
\end{bem}

We can now accurately define the functionals in \eqref{eq:rof-functional-0} and \eqref{eq:rof-functional} for the cases  $\sigma= 0$ or  $\sigma>0$:
\begin{equation}
   \mathcal E(u) = \mathcal E_0(u):=\left\{\begin{array}
    {ll} \displaystyle \mathcal{TV}(u)+\frac{\lambda}{2}\int_\Sigma d_h^2(u,f)\, d\mu_g & {\rm if \ } u\in BV(\Sigma;\mathcal N) \smallskip \\ +\infty & {\rm if\ } u\in L^1(\Sigma;\mathcal N)\setminus BV(\Sigma;\mathcal N)
  \end{array}\right., {\rm and}
\end{equation}

\begin{equation}
  \mathcal E_\sigma(u):=\left\{\begin{array}
    {ll} \displaystyle \mathcal{TV}(u)+\frac{\lambda}{2}\int_\Sigma d_h^2(u,f)\, d\mu_g+\frac{\sigma}{2}\int_\Sigma |du|^2\, d\mu_g & {\rm if \ } u\in H^1(\Sigma;\mathcal N) \smallskip \\ +\infty & {\rm if\ } u\in L^2(\Sigma;\mathcal N)\setminus H^1(\Sigma;\mathcal N)
  \end{array}\right.\,.
\end{equation}
	
\subsection{Euler-Lagrange equations}	

For $\sigma \geq 0$,	a minimizer $u$ of $\mathcal E_\sigma$ has to be a weak solution to the corresponding Euler-Lagrange system:
	\begin{equation}\label{eq:rof-pde}
		\left\{ \begin{array}{rcll}
			\dv_g Z_u & = & -\lambda\exp^{-1}_u\!f &\qquad \text{in $\Sigma$}, \medskip \\
			\nu \cdot Z_u & = & 0 &\qquad \text{on $\partial \Sigma$}, \end{array}\right.
	\end{equation}
	where  $\nu$ represents the outer unit  normal of $\partial\Sigma$, and $Z_u$ is given by
	\[
	Z_u := \bigg(\frac{1}{|d u|} + \sigma\bigg) d u.
	\]
	
Even in the smooth case; i.e $u\in H^1(\Sigma;\mathcal N)$, we need to clarify the meaning of  $\frac{d u}{|d u|}$. Throughout this work, we understand it as a multivalued function defined as
	
	\begin{equation}
		\label{eq:def-Z}
		\frac{d u}{|d u|}: x \longmapsto
		\begin{dcases}
			\frac{d u(x)}{|d u(x)|}, &\text{if $d u(x) \neq 0$,} \smallskip \\
			B_{g\boxtimes h}(0, 1) \subset T_x^\ast \Sigma \otimes T_{u(x)}\mathcal N, &\text{if $d u(x) = 0$.}
		\end{dcases}
	\end{equation}

 Our definition of a (regular) solution to \eqref{eq:rof-pde} is the following one:

	\begin{defn}\label{eq:rof-sol}
		Given $\sigma \geq 0$, we say that $u\in C^{0,1}(\overline\Sigma; \cn)$ is a (regular) solution of \eqref{eq:rof-pde} if there exists $Z\in T^\ast \Sigma \otimes u^\ast T \mathcal N$ with $\dv_g Z \in L^2(\Sigma; \mathbb R^n)$ satisfying
		\begin{equation}
			\label{eq:rof-pde-amb-Z}
			\dv_g Z = -\lambda\exp^{-1}_u \!f \quad \text{with} \quad  Z - \sigma d u\in \frac{d u}{|d u|} \qquad \text{$\mathcal L^2$-a.e. in $\Sigma$},
		\end{equation}
		and it fulfills the homogeneous Neumann condition
		\begin{equation}
			\label{eq:rof-neumann-amb-Z}
			 \nu \cdot Z^\alpha  = 0 \qquad \text{$\mathcal H^1$-a.e. on $\partial\Sigma$},\ \ \text{for all} \quad \alpha=1,\ldots,n.
		\end{equation}
	\end{defn}

	\section{Geodesic convexity of a generalized $p$-energy} \label{convex}
	
	Let us work for this part in a more general framework: $(\mathcal M^{m},g)$ a compact $m$-dimensional manifold with $\partial\mathcal M\neq \emptyset$.	Next, Consider  a  $C^2$ function $F:\mathcal M\times\R\to \R^+$ such that
	\begin{itemize}
		\item $F(x,\cdot)$ is non-decreasing for all $x\in\mathcal M$. \smallskip
		\item $|F(x,s)|\leq C(1+|s|^\frac{p}{2})$ for $(x,s)\in\mathcal M\times \R$ and some $p \geq 1$,  and \smallskip
		\item  the function $G(x,\cdot):= F\left(x,\frac{|\cdot|^2}{2}\right)$ is convex for all $x\in\mathcal M$.
	\end{itemize}
	We consider $\mathcal F: L^2(\mathcal M;\mathcal N) \rightarrow \R$ defined as
	\begin{equation}
	\label{def:functional} \mathcal F(u):=\left\{\begin{array}{ll} \displaystyle \int_{\mathcal M} F\left(\cdot,\tfrac{|du|^2}{2}\right) d\mu_g & {\rm if \ } u\in W^{1,p}(\mathcal M;\mathcal N), \ u(\mathcal M)\subset B_h(p,R) \smallskip \\ +\infty & {\rm otherwise}.\end{array}\right.
	\end{equation}
	In the above setting, the $F$-tension field  is defined by \begin{align}\tau_F(u):&=\sum_i\left(\widetilde \nabla_{e_i}\left(F'\!\left(x,\tfrac{|du|^2}{2}\right)u_*e_i\right)-F'\Big(x,\tfrac{|du|^2}{2}\Big)u_*(\nabla_{e_i}e_i) \right),\end{align}
being $\{e_i\}_{i = 1}^m$ any local orthonormal basis of $T\mathcal M$. Hereafter $F'$ and $F''$ will denote first and second partial derivatives of $F(x,s)$ with respect to $s$.

To prove the convexity of the functional $\mathcal F$, one typically needs to derive formulas for its first and second variation. However, unlike in the previous literature, we have no sort of harmonicity condition which cancels out terms involving the $F$-tension and we have to keep track all extra boundary terms. Taking these difficulties into account, after a tedious but routine computation (which can be reproduced following the lines of \cite{Ara, Ura}), we obtain
	\begin{align}
& \frac{d}{dt}\bigg|_{t=0}\mathcal F(u_t)=-\int_{\mathcal M}\<V,\tau_F(u)\>_h d\mu_g+\int_{\partial\mathcal M}F'\left(x,\tfrac{|du|^2}{2}\right)\<u_\ast\nu,V\>_h\, d\mu_{\tilde g},\end{align}
where $U: \mathcal M \times [0,1] \rightarrow \mathcal N$ is a one-parameter variation of $u \in C^\infty(\mathcal M, \mathcal N)$ so that $U(\cdot, 0) = u$, $u_t :=U(\cdot, t)$ and $V =  U_\ast \partial_t\big|_{t = 0}$.
 Moreover,  $\tilde g$ is the induced Riemannian metric on $\partial \mathcal M$. In turn, the second variation formula reads as
\begin{align}
 \frac{d^2}{dt^2}\bigg|_{t=0}\mathcal F(u_t)  =& -\int_{\mathcal M}\big<\widetilde\nabla_{\partial_t}V,\tau_F(u)\big>_h\, d\mu_g+\int_{\mathcal M}F'\left(x,\tfrac{|du|^2}{2}\right)(|\widetilde\nabla V|^2_h-R(V,du))d\mu_g\\&+\int_{\mathcal M}F''\left(x,\tfrac{|du|^2}{2}\right) \big<\widetilde\nabla V,du\big>_h^2\,d\mu_g+\int_{\partial\mathcal M}
F'\left(x,\tfrac{|du|^2}{2}\right)\big<\widetilde\nabla_{\partial_t}V, u_\ast \nu\big>_h\, d\mu_{\tilde g},
\end{align}
with $$R(V,du)=\sum_{i = 1}^m\, ^h \mathcal R(V,u_*e_i,u_*e_i, V)= \sum_{i = 1}^m {\rm Sec}_h(V,u_*e_i)\big|V\wedge u_*e_i\big|^2_h.$$

To prove the convexity of the energy functional we need the following result about density of smooth functions into manifolds with boundary.
	\begin{lemma}\label{lemma:density}
		Let $R< R_\kappa$. Then, the space  $C^\infty(\mathcal M; \overline{B_h(p,R)})$ is dense in $W^{1,q}(\mathcal M;B_h(p,R))$, for any $q\geq 1$.
	\end{lemma}
	\begin{proof}
		
		Recall that we have an isometric embedding $\iota:  \mathcal{N} \hookrightarrow \R^N$ by Nash.
Then for any $R' \in (R, R_\kappa)$ we can modify the metric outside $\iota(B_h(p,  R')) \subset \R^N$ by flattening it out  (see Figure \ref{fig:flatten}) via a bi-Lipschitz ambient diffeomorphism which keeps the Euclidean metric outside a larger ball so that we get a Riemannian metric on all of $\R^n$ for which the squared distance to the origin is smooth and still strictly convex (for a proof, see \cite[Theorem 5]{PiVe}).

\begin{figure}[h]	
	\includegraphics[scale =0.8]{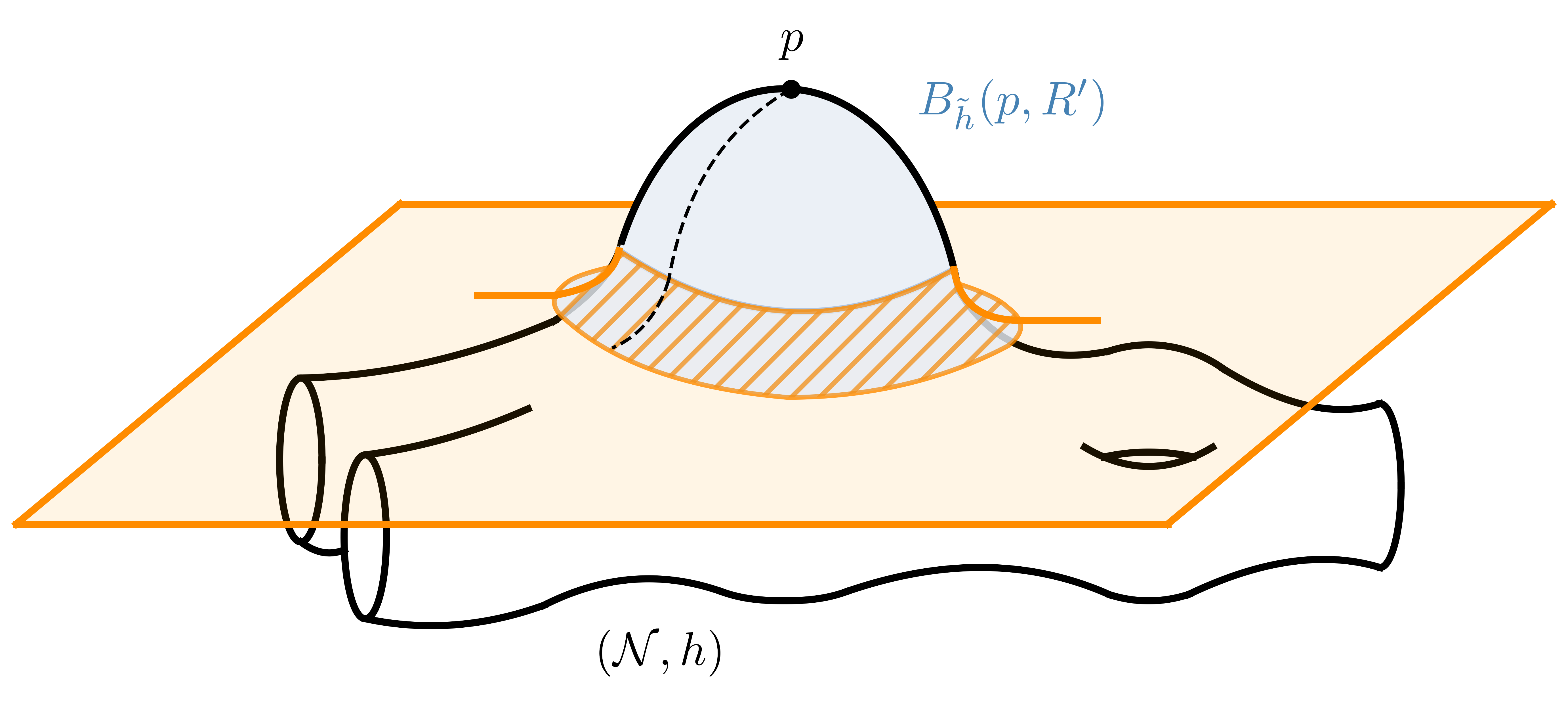}	
	\vspace*{-0.3cm}\caption{Schematic view of the idea of flattening the ends to construct a global coordinate system outside a geodesic ball.} \label{fig:flatten}
\end{figure}


The latter implies that $(\mathcal N, h)$ can be covered by a global, normal coordinate chart centered at $p$, meaning that $(\mathcal N, h)$ can be regarded as $\R^n$ endowed with a metric that can be expressed globally in normal Cartesian coordinates, and hence the extrinsic definition \eqref{Sob-ext} of $W^{1,p}(\mathcal M, \mathcal N)$ coincides with the intrinsic one $W^{1,p}_{\rm in}(\mathcal M, \mathcal N)$ given by \eqref{Sob-in}, cf. \cite[Theorem 4]{PiVe}.
		
		In this setting, given $u \in W^{1,q}(\mathcal M; \mathcal N)$ with $u(\mathcal M) \subset B_h(p,R)$, by \cite[Corollary 7]{PiVe}, we can find a sequence  $u_k \in C^\infty(\mathcal M; \mathcal N)$ such that $u_k \to u$ in $W^{1,q}(\mathcal M; \mathcal N)$ as $k \to \infty$. Accordingly, up to extracting a subsequence, we can assume that the convergence is pointwise
		a.e. and hence that $u_k(\mathcal M) \subset \overline{B_h(p, R)}$.
	\end{proof}

\begin{thm}\label{convexity-p}
	If $\kappa \leq 0$ and $R<R_\kappa$, the functional $\mathcal F$ is geodesically convex.
\end{thm}
\begin{proof}
	We split the proof into several steps:
	
	\emph{ Step 1. Second variation estimate.} In the case that $\kappa\leq 0$ and $U$ is a variation along a geodesic, then
	\[
	 \left.\frac{d^2}{dt^2}\right|_{t=0}\mathcal F({u_t}) \geq \int_{\mathcal M}\left(F'\left(x
	,\tfrac{|du|^2}{2}\right)|\widetilde \nabla V|_h^2+F''\left(x,\tfrac{|du|^2}{2}\right) \big<\widetilde\nabla V,du)\big>_h^2\right)\, d\mu_g\geq 0,\] since $G$ is convex in the second argument.

\emph{ Step 2. Convexity for smooth functions.} Given $u,v\in C^{ \infty}(\mathcal M; \mathcal N)$ with $u(\mathcal M), v(\mathcal M)\subset B_h(p,R)$ for $R<R_\kappa$,  one can always construct (see \cite[proof of Theorem 2.9]{Schoen-an}) a geodesic homotopy $U: \mathcal M \times [0,1] \rightarrow  B_h(p,R)$  with $U(\cdot, 0)=u$ and $U(\cdot, 1)=v$. Let us denote for brevity $U(t) = U(\cdot, t)$ and recall that $t \mapsto U(t)$ are geodesics in $\mathcal N$ for each fixed point in $\mathcal M$. Then, by the previous computations, since $t \longmapsto \mathcal F(U(t))$ is a convex function, it holds:
\begin{equation}\label{convex:smooth}\mathcal F(U(t))\leq t\mathcal F(u)+(1-t)\mathcal F(v)\, \quad {\rm for \ all \ } t\in [0,1].\end{equation}

\emph{Step 3. Convexity for $W^{1,q}(\mathcal M;\mathcal N)$ functions.}
We now observe that, since $B_h(p,R)$ is an NPC-space for $R<R_\kappa$, then $L^2(\mathcal M,B_h(p,R))$ is also an NPC-space \cite[\S 4.1.1]{GigliNobili2021}. Therefore, given $u,v\in L^2(\mathcal M; B_h(p,R))$ there exists a unique geodesic $U(t)$ such that $U(0)=u$, $U(1)=v$. In order to prove that \eqref{convex:smooth} holds  for $u,v$, it suffices to suppose that both belong to $W^{1,q}(\mathcal M;\mathcal N)$ since, otherwise, the right hand side is infinite.

Accordingly, we can apply Lemma \ref{lemma:density} to ensure that we can take $u_k,v_k\in C^\infty(\mathcal M;\overline{B_h(p,R)})$ such that $u_k,v_k\to u,v$ strongly in $W^{1,q}(\mathcal M;B_h(p,R))$. Consider now $U_k$ the geodesic such that $U_k(0)=u_k$ and $U_k(1)=v_k$. Then, since we are in an NPC space we have that \begin{align*}\int_{\mathcal M} \!d_h(U(t),U_k(t))^2\, &\leq (1-t)\int_{\mathcal M}\!d_h(U(t),u_k)^2\, +t\int_{\mathcal M}\!d_h(U(t),v_k)^2 -t(1-t)\int_{\mathcal M}\!d_h(u_k,v_k)^2\, \\& \!\!\stackrel{k\to\infty}\to (1-t)\int_{\mathcal M}d_h(U(t),u)^2+ t \int_{\mathcal M}d_h(U(t),v)^2\,  -t(1-t)\int_{\mathcal M}d_h(u,v)^2\, =0,
\end{align*}
where all the integrals above are computed with respect to $d\mu_g$.

Therefore, by lower semicontinuity with respect to the $L^2$-convergence of the energy and strong convergence of the approximations, we can pass to the limit in \eqref{convex:smooth} and show that \eqref{convex:smooth} is satisfied by $u,v\in W^{1,q}(\mathcal M;B_h(p,R))$.
\end{proof}

This result permits to show that the convexity also holds for the functional $\mathcal{TV}$ given in \eqref{TV-def}.
\begin{thm}\label{convex-tv}
  If $\kappa \leq 0$ and $R<R_\kappa$, then the functional $\mathcal{TV}$ is geodesically convex with $\mathcal N$ replaced by $B_h(p,R)$.
\end{thm}
\begin{proof}
  We let $$\mathcal F_\varepsilon(u):=\left\{\begin{array}{cc}\displaystyle \int_\Sigma \sqrt{\varepsilon+|du|^2}\, d\mu_g & {\rm if \ } u\in W^{1,1}(\Sigma;B_h(p,R))
  \\ +\infty & {\rm otherwise.}\end{array}\right.$$ Then, by Sobolev embedding and Theorem \ref{convexity-p}, we have that $\mathcal F_\varepsilon$ is geodesically convex in $L^2(\Sigma;\mathcal N)$. Note, however, that due to boundedness of both domain and target in the case when the functional is finite, $L^1$ and $L^2$ topologies are equivalent.

  We reason now as in the last part of the proof of Theorem \ref{convexity-p}. We take $u,v\in BV(\Sigma;B_h(p,R))$, then by Theorem \ref{tv-relax}, there exist $\{u_k,v_k\}_k\subset C^1(\Sigma;\overline{B_h(p,R)})$ such that $u_k,v_k\to u,v$ in $L^1(\Omega)$ and $\mathcal{TV}(u_k)\to\mathcal{TV}(u)$ as $k\to\infty$.

   As before, consider $U_k$ the geodesic in $L^2(\Sigma;B_h(p,R))$ such that $U_k(0)=u_k$ and $U_k(1)=v_k$. Then, again since $L^2(\Sigma;B_h(p,R))$ is an NPC space, we obtain that $U_k\to U$ in the $L^2$-topology, with $U$ being the geodesic joining $u$ and $v$. Therefore, letting $\varepsilon=\varepsilon_k\to 0$ as $k\to\infty$, by the lower semicontinuity given by Theorem \ref{tv-relax} and by convexity of $\mathcal F_\varepsilon$, we have
   \begin{align*}\mathcal TV(U)\leq \liminf_{k\to\infty}\mathcal{TV}(U_k(t))& \leq \limsup_{k\to\infty}\left(\mathcal F_{\varepsilon_k}(U_k(t))\right)\leq \limsup_{k\to\infty}\left(t\mathcal F_{\varepsilon_k}(u_k)+(1-t)F_{\varepsilon_k}(v_k)\right) \\ & = t \mathcal{TV}(u)+(1-t)\mathcal{TV}(v). \end{align*}
\end{proof}

	\section{Existence and uniqueness of minimizers} \label{exist}
	
We begin this section with the following invariance principle for minimizers, which guarantees that the small range condition of the source is inherited by any minimizer.

	\begin{prop}\label{prop:comparison}
		Given $\sigma\geq 0$, let $u$ be a minimizer of $\mathcal E_\sigma$. If $f$ satisfies $f(\overline\Sigma) \subset B_h(p, R)$ for some $p\in\mathcal{N}$ and $R < R_\kappa$, then $u(\overline\Sigma)\subset B_h(p, R)$.
	\end{prop}
	
	\begin{proof}
		By assumption, we can introduce geodesic polar coordinates $(r,\theta)$ in $B_h(p,2R)$ centered at $p$. We consider the following Lipschitz retraction: $\pi:\mathcal N\to B_h(p,R)$, with
		\[
		\pi(q):=\begin{cases}
			(r, \theta) & {\rm if  \ } q=(r,\theta)\,,\  r< R\\
			(2R-r, \theta) &  {\rm if  \ } q=(r,\theta)\,,    R \leq r < 2R\\
			p & \text{if } q\in\mathcal N\setminus B_h(p,2R)
		\end{cases},
		\]

We will show that $\pi$ is length decreasing; i.e $d_h(\pi(q),\pi(\tilde q))\leq d_h(q, \tilde q)$ for any $q, \tilde q\in\mathcal N$ as in the proof of \cite[Lemma 10.2.4]{Jost2017}.
First, observe that it is enough to prove
\begin{equation}
  \label{retraction} d_h(\pi(r,\theta),\pi(\tilde r,\tilde \theta))\leq d_h((r,\theta),(\tilde r, \tilde \theta)) \qquad \text{for } \quad q, \tilde q\in B_h(p,2R).
\end{equation}
By construction,  $\pi$ is clearly length decreasing in $r$; i.e. \eqref{retraction} holds with $\theta$ fixed. Therefore, we will only prove the case that $r= \tilde r$, for which
we take a curve $\gamma(s):=(r,\theta(s))$. Now, for each $s$ fixed, $c_s(t):=(t,\theta(s))$ is a radial geodesic with $c_s(0)=p$ and $c_s(r)=\gamma(s)$. Accordingly,
$J_s(t):=\frac{\partial}{\partial s}c_s(t)$ is a Jacobi field with $J_s(0)=0$, $J_s(r)={\gamma}'(s)$ and $d\pi({\gamma}'(s))=J_s(r_\pi)$, with $(r_\pi,\theta)=\pi(\gamma(s))$ ($r_\pi=2R-r<R<r\leq 2R$).

Thus, assuming without loss of generality that $\gamma'(s)\neq 0$, by Rauch comparison theorem and a straightforward computation, we obtain that
\[\frac{|\gamma'(s)|}{|d\pi({\gamma'}(s))|} = \frac{|J_s(r)|}{|J_s(r_\pi)|}\geq \frac{\sin_\kappa(r)}{\sin_\kappa(2R-r)}>1,\]
where $\sin_\kappa(\cdot)=\sin(\sqrt{\max\{\kappa,0\}}\, \cdot)$. Hence $\pi$ is also length decreasing in the $\theta$ directions.

\smallskip
Consequently, $\mathcal E_\sigma(\pi\circ u)\leq \mathcal E_\sigma(u)$, with strict inequality unless $\pi\circ u= u$. Thus, by the minimality of $u$, we conclude that $u(\Sigma) \subset B_h(p, R)$.
	\end{proof}

As a by-product, under the hypothesis of the above statement, we can restrict the search of admissible minimizers of $\mathcal E_\sigma$ to the following set:
 \[\Lambda_R^\sigma :=\left\{\begin{array}{cc} \{u\in H^1(\Sigma;\mathcal N) : u(\Sigma)\subset B_h(p, R)\} & {\rm if \ } \sigma>0 \smallskip \\ \{u\in BV(\Sigma;\mathcal N) : u(\Sigma)\subset B_h(p, R)\} & {\rm if \ } \sigma=0\end{array}\right.\]

		\begin{prop}
			\label{prop:existence}
			Given $\sigma \geq 0$, if $f$ satisfies $f(\Sigma) \subset B_h(p, R)$ for some $p\in\mathcal{N}$ and $R < R_\kappa$, then there exists a minimizer of $\mathcal E_\sigma$ in $\Lambda_R^\sigma$.
	\end{prop}
	\begin{proof}
		We follow the direct method in the calculus of variations. However, we will only give the proof in the case $\sigma=0$, the other case being similar, but easier.

We note that $\mathcal E(u)\geq 0$. Thus there exists a minimizing sequence $\{u_k\}_k$  such that $\mathcal E(u_k)  \to \inf\mathcal E$ as $k \to \infty$. Moreover, by Proposition \ref{prop:comparison}, we can assume that $u_k(\Sigma)\subseteq B_h(p,R)$. Now, by \eqref{boundforiso} and the definition of $\mathcal TV$, we have that $\mathcal E$ is coercive. Therefore, there is a subsequence, not relabeled, and $u\in BV(\Sigma)$ such that $u_k\to u\in L^1(\Sigma)$. Hence we may further assume that $u_k\to u$ a.e. in $\Sigma$. This yields that $u\in\Lambda_R^0$. Finally, by
lower semicontinuity, we conclude that $u$ is a minimizer of $\mathcal E$.
	\end{proof}

As a direct consequence of the convexity results in Theorems \ref{convexity-p}, \ref{convex-tv} and the fact that $d^2_h$ is strictly convex in $B_h(p,R)$ for $R<R_\kappa$, we obtain the following result.

		\begin{prop}
			Under the hypothesis of Proposition \ref{prop:existence}, if $\kappa \leq 0$, then $\mathcal E_\sigma$ has a unique minimizer in $\Lambda_R^\sigma$.
	\end{prop}

	\section{Regularity results for an approximate problem} \label{Hol-reg}
	
	\subsection{Setup of the regularized functional/system} \label{setupA}
	In this section, we will study the regularity of the minimizers of $\mathcal E_\sigma$. Firstly, we  consider the following approximation to $\mathcal E_\sigma$:
	\begin{equation}\label{eq:rof-relaxed-functional}
		\mathcal E_{\sigma, \varepsilon}(u):=\left\{\begin{array}{ll} \displaystyle\int_\Sigma \sqrt{|d u|^2 + \varepsilon^2}\, d\mu_g + \frac{\lambda}{2}\int_\Sigma d_h^2(u,f)\,d\mu_g + \frac{\sigma}{2}\int_\Sigma |d u|^2\, d\mu_g & {\rm if \ } u\in H^1(\Sigma;\mathcal N) \vspace{0.2cm}\\ +\infty & {\rm otherwise.}\end{array}\right.
	\end{equation}
	with $\varepsilon,\sigma>0$. As in Proposition \ref{prop:existence}, one can guarantee the existence of a minimizer for $\mathcal E_{\sigma, \varepsilon}$.
\begin{lemma}
  Given $f\in L^2(\Sigma;\mathcal N)$, there exists $u\in H^1(\Sigma;\mathcal N)$ which minimizes $\mathcal E_{\sigma,\varepsilon}$.
\end{lemma}
\begin{proof}
  Let $\{u_k\}_k\subset H^1(\Sigma;\mathcal N)$ be a minimizing sequence. Then, identifying again $\mathcal N$-valued functions with their extension by Nash embedding, we obtain a subsequence, not relabeled, and a function $u\in H^1(\Sigma)$ such that $u_k\rightharpoonup u $ in $H^1(\Sigma)$ and $u_k\to u$ a.e. in $\Sigma$. Therefore, by lower semicontinuity, we conclude that $u\in H^1(\Sigma;\mathcal N)$ is a minimizer.
\end{proof}

Any minimizer of $\mathcal E_{\sigma,\varepsilon}$ is a weak solution to  the corresponding Euler-Lagrange system
	\begin{equation}\label{eq:relaxed-pde}
		\left\{ \begin{array}{rcll}
			\dv_g \mathcal Z_{\varepsilon, \sigma} & = & -\lambda\exp^{-1}_u f &\qquad \text{in $\Sigma$} \medskip \\
			\nu\cdot  du& = & 0 &\qquad \text{on $\partial\Sigma$}
		\end{array}\right.\tag{$\mathcal S_{\varepsilon, \sigma}^f$}
	\end{equation}
	such that $\mathcal Z_{\varepsilon, \sigma}$ is defined as
	\begin{equation}\label{eq:relaxed-Z}
		\mathcal Z_{\varepsilon, \sigma} =\mathcal Z_{\varepsilon, \sigma}(x,du) := \bigg(\frac{1}{\sqrt{|du|^2 + \varepsilon^2}} + \sigma\bigg) du.
	\end{equation}

Now the goal is to study the regularity of weak solutions to \eqref{eq:relaxed-pde}. In order to do this, we write the system locally using extrinsic coordinates in the target, by Nash theorem and isothermal coordinates in the domain, as described in Sections \ref{tensions} and \ref{p-energy-sect}. Moreover, since we want to obtain regularity estimates up to the boundary of $\Sigma$, we will focus on neighborhoods around any point in $\partial\Sigma$, the proofs for interior neighborhoods being similar, but easier than the boundary case.

We now fix some notation. For $x\in\Sigma$ and $r > 0$, $B_r(x)$ denotes the open ball in $\R^2$ with center $x =(x_1, x_2)$ and radius $r >0$.  We consider the following sets:
	\begin{align*}
		B_r^+(\overline x) &= B_r( x)\cap \{x_2 >  \overline x_2\}\\
		\Omega_r(\overline x) &= B_r(\overline x) \cap B_1^+(0)
	\end{align*}

The first step towards regularity is to flatten $\partial\Sigma$ around a
fixed point $z\in\partial \Sigma$. As $\partial \Sigma$ is smooth, there exists a smooth transformation $\Psi_{z}$ for a some neighborhood $\mathcal U$ of $z$ such that
		\begin{equation*}
			\Psi_{z}(\mathcal U\cap\Sigma) = B_1(0) \cap\{x_2 > 0\} =: B_1^+, \quad \Psi_{z}(\mathcal U\cap \partial \Sigma) = B_1(0) \cap\{x_2 = 0\} =: \Gamma_1.
		\end{equation*}
\begin{figure}[h]	
			\includegraphics[scale =0.9]{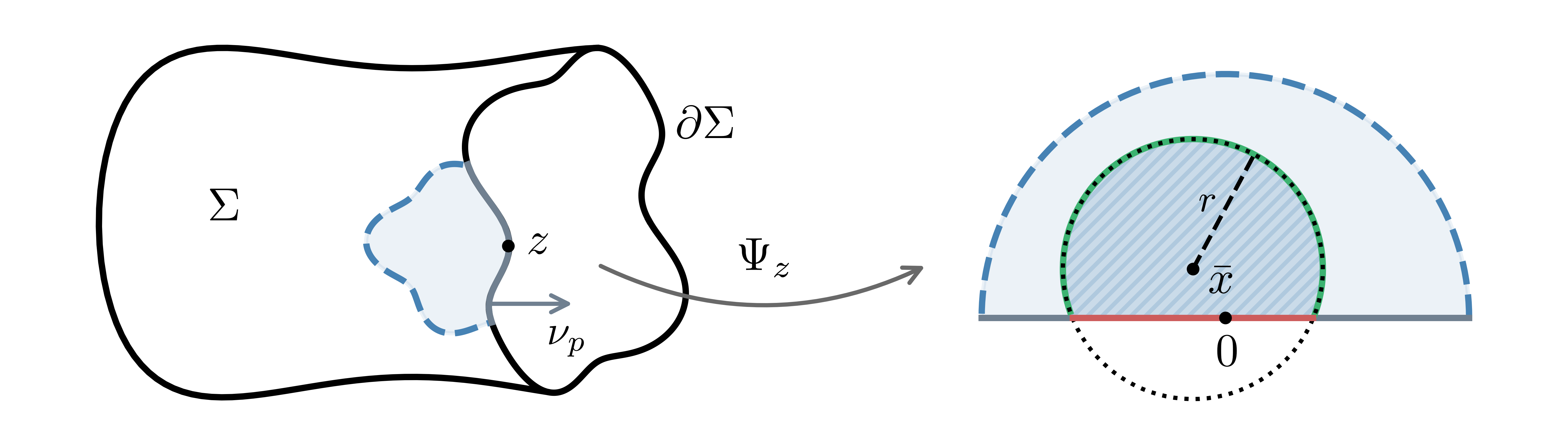}	
			\vspace*{-0.3cm}\caption{Setup for the problem after locally flattening the boundary. The blue area denotes the neighborhood $\mathcal{U}$ in the left figure and $B_1^+$ in the right figure. The striped area represents $\Omega_r(\bar{x})$, the red segment indicates $\Gamma_1\cap \partial\Omega_r$, and the green segment marks the region where $u=v$.}
		\end{figure}

Without loss of generality, we may
assume that $\mathcal U= U_\ell$ and $\Psi_z=\phi_\ell$ for some $\ell$. For simplicity we write $\varrho:=\varrho_\ell$ (recall the notation from Section \ref{p-energy-sect}). 
Therefore, any weak solution to $(\mathcal S_{\varepsilon,\sigma}^f)$ must satisfy \begin{equation}
			\label{eq:simplified-pde}
			\begin{dcases}
				\begin{aligned}
					\sum_{i}\partial_i a_i^\alpha(x,\nabla u) &= \Phi^\alpha(x, u,\nabla u) \quad \text{for $x\in B_1^+$}\\
					a^\alpha_2(x,\nabla u)|_{\Gamma_1}&= 0
				\end{aligned}
			\end{dcases}
		\end{equation}
		for $\alpha\in \{1, ..., N\}$, %
with
$a_i^\alpha(x,\xi)=(b(x,\xi)+\sigma)\xi_i^\alpha$, whereas $b$ and $\Phi$ are defined as
	\begin{align*}
		b(x,\xi) &:= \frac{\varrho(x)}{\sqrt{|\xi|^2 + \varepsilon^2\varrho^2(x)}} ,\\
		\Phi(\cdot,p,\xi) &:= -(b(\cdot,\xi)+\sigma)\sum_{i}\mathcal{A}_{p}(\xi_i, \xi_i)-\lambda\varrho^2\exp^{-1}_{p} f.
	\end{align*}

In particular, we note that $a$ satisfies by definition that
	
	{\bf {\small (H1)}}	$\displaystyle |a(x,\xi)| \leq \left(\frac{C}{\varepsilon} + \sigma\right)|\xi|$.
	
	Moreover, the matrix $A^{\alpha\beta}_{i j}(x,\xi) := \frac{\partial a^\alpha_i(x,\xi)}{\partial \xi^\beta_j}$ satisfies the following ellipticity condition:
	\label{eq:ellipticity}
	
	{\bf {\small (H2)}}	$\displaystyle \sum_{i,j,\alpha,\beta}A^{\alpha \beta}_{ij}(x,\xi)\eta^\alpha_i\eta^\beta_j  =  \frac{-b(x,\xi)}{\xi^2 + \varrho^2(x)\varepsilon^2}\sum_{i,\alpha}(\xi^\alpha_i\eta^\alpha_i)^2 + (b(x,\xi) + \sigma)|\xi|^2 \geq \sigma|\eta|^2$
for all $\eta\in\mathbb R^{2N}$.

Furthermore, for every pair $p,q\in\mathbb R^{2N}$ the next inequalities hold:
	\begin{enumerate}
		\item[{\bf {\small (H3)}}]	$\displaystyle \sum_{i,\alpha}(a_i^\alpha(x,\xi) - a_i^\alpha(x,\eta))(\xi^\alpha_i - \eta^\alpha_i)  \geq \sigma |\xi - \eta|^2,$

		\item[{\bf {\small (H4)}}]	$\displaystyle	\sum_{i,\alpha} a_i^\alpha(x,\xi)\xi^\alpha_i \geq \sigma |\xi|^2.$

 		\item[{\bf {\small (H5)}}] 	$\displaystyle |a(x,\xi)-a(y,\xi)| \leq \frac{C}{\min \varrho}|x-y|$
		\item[{\bf {\small (H6)}}] 	$\displaystyle |a(x,\xi)-a(x,\eta)|  \leq \left(\frac{C}{\min \varrho}+\sigma\right)|\xi-\eta|$.
	\end{enumerate}

Now, we observe that $$\Phi(x,p,\xi)={\bf n}(x,p,\xi)+{\bf t}(x,p)\in T_p\mathcal N\oplus T_p^\bot\mathcal N,$$
with ${\bf t}(x,p)=-\lambda\varrho^2\exp^{-1}_{p} f.$

Therefore, we can unconstrain the problem by extending ${\bf n}$ and ${\bf t}$ to $B_1^+\times\R^N\times\R^{2N}$ as in \cite[Section 2]{CRMP}. Thus, we assume that ${\bf n}$, ${\bf t}$ (and therefore $\Phi$) have already been extended in such a way that the following holds:

\begin{enumerate}
\item[{\bf {\small (H7)}}] 	$|{\bf t}(x,p)|\leq C|p|$, $|\Phi(x,p,\xi)|\leq C_1+C_2|\xi|^2$.
\end{enumerate}
	
\medskip
We will also need a generalized Poincar\'e inequality for hemispheres. The proof of the following result, which we omit, is an easy adaptation of the Poincar\'e inequality for balls in \cite[Theorem 2 p.291]{Evans}, combined with Lemma \ref{lem:poincare}.
\begin{lemma}\label{lem:poincareball}Let $\overline x\in B_1^+$ and $r_0<{\rm dist}(\overline x,\partial B_1^+\setminus \Gamma_1)$. Then, there exists $C>0$ such that, for all $r\leq r_0$, it holds \begin{itemize}
\item [(a)] $\displaystyle
\|u\|_{L^1(\Omega_r)}\leq C r\|\nabla u\|_{L^1(\Omega_r)},\quad {\rm for \ any \ } u\in W^{1,1}_{\partial \Omega_r\setminus\Gamma_1}(\Omega_r)$,
\item [(b)] $\displaystyle
\|u-(u)_{\overline x,r}\|_{L^1(\Omega_r)}\leq C r\|\nabla u\|_{L^1(\Omega_r)},\quad {\rm for \ any \ } u\in W^{1,1}(\Omega_r)$
\end{itemize}
\end{lemma}

\medskip	
	\subsection{Higher integrability of solutions}

Let $\overline x\in B_1^+$ and $r<{\rm dist}(\overline x,\partial B_1^+\setminus \Gamma_1)$ and suppose that $u(\overline \Omega_r(\overline x))\subset \overline{B_g(p, R)}$ with $R<R_\kappa$. We define the barycenter of $u$ in  $\Omega_r(\overline x)$ as the barycenter of the measure $\mu:=u\#\mathcal L^m|_{\overline{B_r(x_0)}},$; that is, it can be regarded just as the minimizer of
\[\mathfrak b \mapsto \int_{\Omega_r(\overline x)} d^2_{\vp h}(u(x), \mathfrak b)\, dx.\]

Let
\begin{equation} \label{Rast}
R_\kappa^*:=\left\{\begin{array}
{cc} \frac{{\rm inj}_p\mathcal N}{2} & {\rm if \ } \kappa\leq 0 \medskip \\ \min\left\{\frac{\mathrm{inj}_p \mathcal{N}}{2}, \frac{\pi}{4 \sqrt{\kappa}}\right\} & {\rm if \ } \kappa>0
\end{array}\right..
\end{equation}
With the same proof as \cite[Propositions 1 and 2]{CRMP} we obtain the following result (we omit the details).

	\begin{prop}\label{prop:CMP-inequalities}
		Let $f\in L^2(\Sigma;\mathcal N)$ be such that $f(\Sigma)\subset B_h(p,R)$ with some $p\in\mathcal N$ and $R<R_\kappa^*$ and let $u\in\Lambda_R^\sigma$ be a weak solution of \eqref{eq:relaxed-pde}. Then the next intrinsic Caccioppoli inequality is satisfied
		\begin{equation}\label{eq:intrinsic-Cacciopoli}
			\int_{\Omega_r(\overline x)}|\nabla u|^2dx \leq \frac{C}{r^2}\int_{\Omega_{2r}(\overline x)}(1 + \dg^2(u, \mathfrak{b}_{2 r}))dx
		\end{equation}
		for all $\overline x\in\overline B_1^+$ and $r < \frac{\text{dist}(\overline x, \partial B_1^+\setminus \Gamma_1)}{2}$, where $C > 0$ and $\mathfrak{b}_{2 r}$ denotes the barycenter of $u$ in $\Omega_{2r}(\overline x)$. Consequently, there exist $q > 2$ and $\tilde C > 0$ such that $u\in W^{1,q}(B_1^+)$ and
		\begin{equation}\label{eq:sobolev-bound}
			\left(\fint_{\Omega_r(\overline x)}(1 + |\nabla u|)^qdx\right)^{1/q} \leq \tilde C \left(\fint_{\Omega_{2r}(\overline x)}(1 + |\nabla u|)^2dx\right)^{1/2}
		\end{equation}
		for all $\overline x\in\overline B_1^+$ and $r < \frac{\text{dist}(\overline x, \partial B_1^+\setminus \Gamma_1)}{2}$.
	\end{prop}

From now on, we will assume that $f(\Sigma)\subseteq B_h(p,R)$ for some $p\in\mathcal N$ and $R<R_\kappa^*$.

	\subsection{Hölder regularity of the approximate solutions}
	
	The next step is to prove H\"older regularity of weak solutions to \eqref{eq:relaxed-pde}.

	\begin{thm}\label{thm:holder-u}
		Let $u\in H^1(\Sigma;\mathcal N)$ be a weak solution of \eqref{eq:relaxed-pde}. Then $u\in C^{0,\beta}(\overline \Sigma; \mathcal N)$ for all $\beta\in(0,1)$.
	\end{thm}
	\begin{proof}

Fix $\overline x\in \overline{ B_{1-\delta}^+(0)}$, with $0<\delta<1$, take $r < \frac{\delta}{2}<1$, and consider $v$ a solution to the system
		\begin{equation}\label{eq:hom-pde}
			\begin{dcases}
				\begin{aligned}
					\sum_{i}\partial_i a_i^\alpha (\overline x,\nabla v) = 0, & \quad \text{in $\Omega_r(\overline x)$}, \quad \alpha =1, \ldots N,\\
					a_2^\alpha(\overline x,\nabla v)|_{\partial\Omega_r(\overline x)\cap \Gamma_1}= 0 &, \quad v|_{\partial\Omega_r(\overline x)\setminus \Gamma_1} = u.
				\end{aligned}
			\end{dcases}
		\end{equation}
		We note that it is easy to adapt \cite[Theorem 3.I]{Campanato1987} to our case, which leads to the inequality
		\begin{equation}\label{eq:campanato-estimate-1}
			\int_{\Omega_\rho(\overline x)}|\nabla v|^2dx \leq C\left(\frac{\rho}{r}\right)^2\int_{\Omega_r(\overline x)}|\nabla v|^2dx
		\end{equation}
		for all $\rho \leq r$ and some $C > 0$. Hereinafter, we will omit $\overline x$ if it is not needed. As $v$ is a weak solution of \eqref{eq:hom-pde},  it satisfies
		\begin{equation*}
			\sum_{i,\alpha}\int_{\Omega_r} a^\alpha_i(\overline x,\nabla v) \partial_i \varphi^\alpha \, dx = 0
		\end{equation*}
		for all $\varphi\in {H}^1_{\partial\Omega_r\setminus\Gamma_1}(\Omega_r)\cap L^\infty(\Omega_r)$. If $\varphi = v - u$, by (H1) and (H4) one can estimate
		\[\sigma \int_{\Omega_r} |\nabla v|^2 dx \leq \sum_{i,\alpha}\int_{\Omega_r} a_i^\alpha (\overline x, \nabla v) \partial_i u^\alpha \, dx \leq \big(\frac{C}{\varepsilon} + \sigma\big) \int_{\Omega_r} |\nabla v||\nabla u|\, dx,\]
		which by means of Young's inequality yields
		\begin{equation}\label{eq:DvDu}
			\int_{\Omega_r} |\nabla v|^2dx \leq C\int_{\Omega_r} |\nabla u|^2dx
		\end{equation}
		for some $C > 0$ which depends on $\sigma$ and $\varepsilon$.
		
		On the other hand, from equations \eqref{eq:simplified-pde} and \eqref{eq:hom-pde}, we deduce that $u$ and $v$ fulfill
		\begin{align}
			\label{eq:arkhi-weak-sol}
			\sum_{i,\alpha}\int_{\Omega_r}(a_i^\alpha(\overline x,\nabla u) - a_i^\alpha(\overline x,\nabla v)) \partial_i \varphi^\alpha \,  dx  = &
			\sum_{\alpha}\int_{\Omega_r} \Phi^\alpha(x, u, \nabla u) \varphi^\alpha dx \\ & +\sum_{i,\alpha} \int_{\Omega_r} (a^i_\alpha(\overline x,\nabla u) - a^i_\alpha(x,\nabla u))\partial_i \varphi^\alpha \, dx \nonumber
		\end{align}
		for all $\varphi\in {H}^1_{\partial\Omega_r\setminus\Gamma_1}(\Omega_r)\cap L^\infty(\Omega_r)$. In particular, set $w = u - v$ and take $M > 0$ so that
		\begin{equation} \label{m-bdd}
		|w|_{L^\infty(\Omega_{r/2})} \leq M.
		\end{equation}
		Let us consider $\varphi = w\cdot(T^s - (|w|+M)^s)_+$  with $s\in(0,1)$ and $T > 0$ to be specified later on. Then, by (H4) and (H6) the left side of \eqref{eq:arkhi-weak-sol} can be estimated as follows
		\begin{align*}
			\sum_{i,\alpha}\int_{\Omega_r} (a_i^\alpha(\overline x,\nabla u) - a_i^\alpha(\overline x, \nabla v)) \partial_i\varphi^\alpha\, dx \geq &
			\, \sigma \int_{\Omega_r} |\nabla w|^2(T^s - (|w|+M)^s)_+dx \\ &- \big(\frac{C}{\min \varrho} + \sigma\big)s\int_{K_T} |\nabla w|^2 (|w|+M)^{s}dx
		\end{align*}
		where $K_T := \{x \in \Omega_r : |w(x)| < T - M\}$. Concerning the right side of \eqref{eq:arkhi-weak-sol}, from (H7)  we deduce that
		\begin{equation}
			\label{eq:arkhi-starting-ineq-thm1}
			\sum_{\alpha}\int_{\Omega_r}|\Phi^\alpha(x,u,\nabla u)\varphi^\alpha|dx \leq \int_{\Omega_r}(C_1 + C_2|\nabla u|^2)|w|(T^s - (|w| + M)^s)_+dx,
		\end{equation}
		while (H5) implies
		\begin{align}
			\label{eq:arkhi-starting-ineq-thm1.1}
			\sum_{i,\alpha}\int_{\Omega_r} (a_i^\alpha(\overline x,\nabla u) - a_i^\alpha(x,\nabla u)) \partial_i\varphi^\alpha\, & dx\  \leq \frac{C}{\min\varrho}\int_{\Omega_r}|x-\overline x||\nabla\varphi|\, dx\\
			 & \leq  C\int_{\Omega_r}|x-\overline x||\nabla w|(T^s-(|w|+M)^s)_+ dx \medskip \\ & \quad +Cs\int_{K_T}|x-\overline x||\nabla w|(|w|+M)^s dx.
		\end{align}
	
		Now, plugging the above estimates into \eqref{eq:arkhi-weak-sol}, taking into account that $|w|+ M < T$ on $K_T$ and $|x -\bar x| < r$, we derive
		\begin{align*}
			\sigma \int_{\Omega_r} |\nabla w|^2(T^s & - (|w|+M)^s)_+dx
			\leq (C + \sigma)sT^s\int_{K_T}|\nabla w|^2 dx \\ & + \int_{\Omega_r}(C_1 + C_2|\nabla u|^2)|w|(T^s - (|w| + M)^s)_+dx\\ & +C r \int_{\Omega_r}|\nabla w|(T^s-(|w|+M)^s)_+ dx +C r s T^s\int_{K_T}|\nabla w| dx,
		\end{align*}

		Young's inequality applied to the last two terms, joint with $r, s <1$, yields
		\begin{align*}
			\frac{\sigma}{2} \int_{\Omega_r} |\nabla w|^2(T^s - (|w|+M)^s)_+dx
			&\leq CsT^s\int_{K_T}|\nabla w|^2 dx +CT^sr^2 \\ & + C_1 T^s\int_{\Omega_r} |w| dx + C_2 T^s\int_{K_T} |\nabla u|^2|w|dx.
		\end{align*}
		Now we take $T=2^{\frac{1}{s}+1} \, M$ so that, by definition of $M$ in \eqref{m-bdd}, one gets
		\begin{equation*}
			\sigma \int_{\Omega_r} |\nabla w|^2(T^s - (|w|+M)^s)_+dx \geq \frac{\sigma T^s}{2}\int_{\Omega_{r/2}}|\nabla w|^2dx.
		\end{equation*}
		Then we deduce that
		\begin{equation}\label{eq:Dw-sense-Ts}
			\int_{\Omega_{r/2}} |\nabla w|^2dx \leq Cs\int_{\Omega_r}|\nabla w|^2dx + C_1\int_{\Omega_r}|w|dx + C_2\int_{K_T}|w||\nabla u|^2dx +Cr^2
		\end{equation}
		for some $C > 0$.
		In the sequel we will estimate the integrals of the right side. For the first one, from \eqref{eq:DvDu} we reach
		\begin{equation*}
			\int_{\Omega_r}|\nabla w|^2dx \leq  \int_{\Omega_r} \big(|\nabla u|^2 + |\nabla v|^2\big)dx \leq  C \int_{\Omega_r} |\nabla u|^2 \, dx.
		\end{equation*}
		To estimate the second integral, by means of Poincaré's inequality,  Lemma \ref{lem:poincareball}(a), Hölder and the previous inequality,  we observe that
	\[		\int_{\Omega_r}|w|\,dx \leq Cr\int_{\Omega_r}|\nabla w| dx \leq \widetilde Cr^2\left(\int_{\Omega_r} |\nabla u|^2 \, dx\right)^{1/2} \leq C_R r^2,\]
		where we have also applied the Caccioppoli inequality in \eqref{eq:intrinsic-Cacciopoli}. Finally, for the third integral we have that
		\begin{align*}
			\int_{K_T}|w||\nabla u|^2dx 
			&\leq Cr^2\left(\fint_{\Omega_r}|\nabla u|^q dx\right)^{\frac{2}{q}}\left(\fint_{\Omega_r}|w|^{\frac{q}{q-2}} dx\right)^{\frac{q-2}{q}}
			\\
			& \leq CT^{\frac{2}{q}}r^2\fint_{\Omega_{2r}}(1 + |\nabla u|^2)dx\left(\fint_{\Omega_r}|w| dx \right)^\frac{q-2}{q}  \\ &\leq CT^\frac{2}{q}\left(\int_{\Omega_r}|\nabla u|^2 dx\right)^\frac{q-2}{2q}\int_{\Omega_{2r}}(1+|\nabla u|^2)dx
		\end{align*}
		using H\"older's and \eqref{eq:sobolev-bound} inequalities and the definition of $K_T$. Thus, we deduce from \eqref{eq:Dw-sense-Ts} that
		\begin{equation}\label{eq:Dw-bound-neta}
			\int_{\Omega_{r/2}} |\nabla w|^2dx \leq C\left(\left(s + T^\frac{2}{q}\left(\int_{\Omega_{r}}|\nabla u|^2dx\right)^\frac{q-2}{2q}\right)\int_{\Omega_{2r}}|\nabla u|^2dx + (1+T^\frac{2}{q}){r^2}\right)
		\end{equation}
		If we define $\phi(\rho) := \int_{\Omega_\rho}|\nabla u|^2dx$, we deduce for all $\rho < r/2$ that
		\begin{align*}
			\phi(\rho) & \leq C\left(\int_{\Omega_{r/2}}|\nabla w|^2dx + \left(\frac{\rho}{r}\right)^2\int_{\Omega_{r/2}}|\nabla v|^2dx\right)\\
			&\leq C\int_{\Omega_{r/2}}|\nabla w|^2dx + C\left(\frac{\rho}{r}\right)^2 \phi(r) \\
			&\leq C\left(s + \left(\frac{\rho}{r}\right)^2 + T^\frac{2}{q}\phi(r)^\frac{q-2}{2q}\right)\phi(2r) + C(1+T^\frac{2}{q})r^2
		\end{align*}
		applying \eqref{eq:campanato-estimate-1}, \eqref{eq:DvDu} and \eqref{eq:Dw-bound-neta} inequalities.
		Possibly increasing $C$ we deduce, for $\rho<r/4$,
		\begin{equation*}
			\phi(\rho)
			\leq C\left(s + \left(\frac{\rho}{r}\right)^2 + T^\frac{2}{q}\phi(r)^\frac{q-2}{2q}\right)\phi(r) + C(1+T^\frac{2}{q})r^2.
		\end{equation*}
		Observe that for $r/4<\rho<r$, we obviously have $\phi(\rho)\leq \phi(r)\leq 16\left(\frac{\rho}{r}\right)^2\phi(r)$. Then, the previous inequality holds for any $\rho<r$.
		
		Setting $\rho=\tau r$, $\beta\in (0,1)$ and $s=\tau^2$, with $3C\tau^2= \tau^{2\tilde\beta}$, $\tilde\beta\in (\beta,1)$. Since $T^{\frac{2}{q}}\phi(r)^\frac{q-2}{q}\leq \tau^2$, for $\phi(r)\leq \varepsilon_1$, we get
		\begin{equation}
			\phi(\tau r)
			\leq \tau^{2\tilde\beta}\phi(r) + C(1+T^\frac{2}{q})r^2\leq \tau^{2\tilde\beta}\phi(r) + C(1+T^\frac{2}{q})r^{2\beta}.
		\end{equation}
		Therefore, by induction, we get
		\begin{align}
			\phi(\tau^k r)&\leq \tau^{2\tilde\beta k}\phi(r)+C(1+T^\frac{2}{q})(\tau^{k-1}r)^{2\beta}\sum_{s=0}^\infty (\tau^{2(\tilde\beta-\beta)})^s\\ &=\left(\phi(r)+C(1+T^\frac{2}{q}) \frac{r^{2\beta}}{\tau^{2\beta}-\tau^{2\tilde\beta}}\right)\tau^{2\beta k}<\varepsilon_1\tau^{2\beta k},
		\end{align}
		if we choose $\phi(r)\leq \frac{\varepsilon_1}{2}$ and $C(1+T^\frac{2}{q}) \frac{r^{2\alpha}}{\tau^{2\alpha}-\tau^{2\beta}}\leq \frac{\varepsilon_1}{2}$, which holds for $r\leq r_0$ for some $r_0<\frac{\delta}{2}$. Then,
		\begin{equation}
			\phi(\rho)\leq C\left(\frac{\rho}{r_0}\right)^{2\beta},\quad {\rm for \ every \ }\rho\leq r_0.
		\end{equation}
		
		Applying Poincaré's inequality, Lemma \ref{lem:poincareball}(b) it leads to $u\in \mathcal L^{2+2\beta}(\Omega_r)$, and thus $C^{0,\beta}(\Omega_{r_0})$ for all $\beta\in (0,1)$  thanks to \cite[Theorem I.2]{Campanato1963}.  Therefore, $u \in C^{0,\beta}(\overline\Sigma; \mathcal N)$ for all $\beta\in(0,1)$.
	\end{proof}
	
	\subsection{Hölder regularity of the gradient}
	\begin{thm}
		\label{thm:holder-Du}
		Let $u$ be a weak solution of \eqref{eq:relaxed-pde}. Then $u\in C^{1,\beta_0}(\overline{\Sigma}; \mathcal N)$ for some $\beta_0 \in (0,1)$.
	\end{thm}
	\begin{proof}
		Again, for $0<\delta<1$, we fix $\overline x\in \overline{ B_{1-\delta}^+(0)}$, and take $r < \frac{\delta}{2}$. Recall that Theorem \ref{thm:holder-u} ensures that $u\in C^{0,\beta}(\overline \Sigma; \mathcal N)$ for every  $\beta\in(0,1)$ and thus, $u$ satisfies
		\begin{equation}
			\label{eq:arkhi-r-bounds}
			\osc_{\Omega_r(\overline x)} u \leq C r^{\beta}, \quad \int_{\Omega_r(\overline x)}|\nabla u|^2{\vp dx} \leq C r^{2\beta}.
		\end{equation}
		for $r < \frac{\delta}{2}$ and $\beta\in(0,1)$.
		
		Let $v$ be again a solution to \eqref{eq:DvDu}, and set $w :=  u - v$. Now we will apply \eqref{eq:arkhi-weak-sol} with $\varphi = w\min\{1, \ell/|w|\}$  for $\ell \geq 1$. From the left side of $\eqref{eq:arkhi-weak-sol}$, we deduce that
		\begin{align*}
			\sum_{i,\alpha}\int_{\Omega_r}(a_i^\alpha(\overline x,\nabla u) - a_i^\alpha(\overline x,\nabla v)) \partial_i\varphi^\alpha \, dx 
			\geq & \, \sigma\int_{\Omega_r}|\nabla w|^2 \min\left\{1, \frac{\ell}{|w|}\right\}dx \\& - (C+\sigma) \int_{K_\ell} |\nabla w|^2dx,
		\end{align*}
		where $K_\ell := \{x\in\Omega_r : |w(x)| > \ell\}$.
		
		From the right side of \eqref{eq:arkhi-weak-sol}, by means of (H7), we have that
		\begin{align*}
			\sum_\alpha\int_{\Omega_r} |\Phi^\alpha(x,u,\nabla u) \varphi^\alpha|dx & \,  {\leq} \int_{\Omega_r} (C_1 + C_2|\nabla  u|^2)|\varphi|dx\\
			& \leq  C_1 \int_{\Omega_r}|w|dx + C_2\int_{\Omega_r}|\nabla  u|^2\min\left\{|w|, \ell\}\right\}dx
		\end{align*}
		
		Moreover, using (H5) and Young's inequality, we can write
		\begin{align}\sum_{i,\alpha}\int_{\Omega_r} (a_i^\alpha(\overline x, &\nabla u) - a_i^\alpha(x,\nabla u))\partial_i \varphi^\alpha \, dx\leq C\int_{\Omega_r}|x-\bar x||\nabla \varphi| dx \\ & \leq Cr\left(\int_{\Omega_r}|\nabla w|\min\left\{1,\frac{\ell}{|w|}\right\} dx+\int_{K_\ell}|\nabla w| \,dx\right) \leq \tilde C r \int_{\Omega_r} |\nabla w| \, dx\\& \leq \frac{\sigma}{4} \int_{\Omega_r} |\nabla w|^2 \, dx + C_\sigma r^4.
		\end{align}

	Consequently, plugging the above estimates into \eqref{eq:arkhi-weak-sol}, we reach the next inequality
			\begin{align}   \label{eq:arkhi-1}
				\sigma\int_{\Omega_r}|\nabla w|^2 \min &\left\{1,  \frac{\ell}{|w|}\right\}  dx \leq  C_1\int_{\Omega_r}|w|dx + C_2\int_{\Omega_r}|\nabla u|^2\min\{|w|, \ell\}dx \nonumber \\
				& \hspace*{1.0cm} + C\int_{K_\ell}|\nabla w|^2dx +  \frac{\sigma}{4} \int_{\Omega_r} |\nabla w|^2 \, dx + C_\sigma r^4.
			\end{align}	
		
		Now we note that, applying Lemma \ref{lem:poincareball}(a), Hölder and Young for any $\theta > 0$, it holds
		\begin{align}\label{eq:arkhi-q-bounds}
				\int_{\Omega_r} \!\! |w|dx & \leq  Cr\int_{\Omega_r} |\nabla w|dx 
				\leq \theta\int_{\Omega_r} \!\!|\nabla w|^2dx + C_\theta r^4
		\end{align}
	
	On the other hand, notice that
\[\int_{\Omega_r}|\nabla w|^2 \min \left\{1,  \frac{\ell}{|w|}\right\}  dx \geq \int_{K_\ell^c} |\nabla w|^2 dx = \int_{\Omega_r} |\nabla w|^2 dx - \int_{K_\ell} |\nabla w|^2 dx.\]	
Substituting this and \eqref{eq:arkhi-q-bounds} into \eqref{eq:arkhi-1} and taking $C_1\theta = \frac{\sigma}{4}$, we attain
			\begin{align}
				\label{eq:arkhi-part-1}
				\frac{\sigma}{2}\int_{\Omega_r}|\nabla w|^2dx  \leq (C+\sigma)\int_{K_\ell}|\nabla w|^2dx
				+ C_2\int_{\Omega_r}|\nabla u|^2\min\{|w|, \ell\}dx + C r^4
			\end{align}
		Parallel to this, we recall that $u$ satisfies
		\begin{equation*}
			\sum_{i,\alpha}\int_{\Omega_r}a_i^\alpha(x,\nabla u)\partial_i \varphi^\alpha \, dx = \sum_{\alpha}\int_{\Omega_r} \Phi^\alpha(x,u,\nabla u)\varphi^\alpha(x)dx,
		\end{equation*}
		for $\varphi\in {H}^1(\Omega_r\cup\gamma_r)\cap L^\infty(\Omega_r)$. By setting $\varphi := (u - u(\overline x))\min\{|w|, \ell\}$ and using the notation $\partial_i w = (\partial_i w^1, \ldots, \partial_i w^N)$, we have the equality	
			\begin{align} 	\label{eq:arkhi-u-varphi}
				\sum_{i,\alpha}\int_{\Omega_r} a_i^\alpha(x,\nabla u) \partial_i u^\alpha  \min&\{|w|, \ell\} dx  + \sum_{i,\alpha}\int_{K_\ell^c} a_i^\alpha(x,\nabla u)(u^\alpha - u^\alpha(\overline x)) \frac{1}{|w|} w \cdot \partial_i w\, dx \\ 
				&= \sum_{\alpha}\int_{\Omega_r} \Phi^\alpha(x,u, \nabla u)(u^\alpha - u^\alpha(\overline x))\min\{|w|, \ell\}dx \\
				&\leq \osc_{\Omega_r}u \left(C_1\int_{\Omega_r} |w|dx + C_2\int_{\Omega_r}|\nabla u|^2\min\{|w|, \ell\}dx\right),
			\end{align}
		where the latter follows from (H7). Thanks to (H4) and (H1), applied to the first and second term on the left hand side above, we reach
		\begin{align*}
			\sum_{i,\alpha}\int_{\Omega_r} a_i^\alpha(x,\nabla u) &\partial_i u^\alpha\min\{|w|, \ell\} dx \, {\geq} \,  \sigma\int_{\Omega_r}|\nabla u|^2\min\{|w|, \ell\} dx,\\
			-\sum_{i,\alpha}\int_{K_\ell^c} a_i^\alpha(x,\nabla u) &(u^\alpha - u^\alpha(\overline x)) \frac{w \cdot \partial_i w}{|w|} dx \, {\leq} \, \big(\frac{C}{\varepsilon} + \sigma\big)\osc_{\Omega_r} u \int_{\Omega_r}|\nabla u||\nabla w|dx.\\
		\end{align*}
	
	\vspace*{-0.4cm}
		Then we deduce that there exists a $C > 0$ such that
		\begin{equation}\label{eq:arkhi-previa-epsilon-zeta}
			\begin{aligned}
				(\sigma - C_2 \osc_{\Omega_r}u)&\int_{\Omega_r} \!\!\!|\nabla u|^2\min\{|w|, \ell\} \, dx
				\leq \osc_{\Omega_r} u \bigg(C_1  \!\! \int_{\Omega_r}  \!\!\! |w|dx + C  \int_{\Omega_r}  \!\!\!|\nabla u||\nabla w|dx\bigg)
			\end{aligned}
		\end{equation}
		Applying Young's inequality for any  $\zeta > 0$ and both estimates in \eqref{eq:arkhi-r-bounds}, we find for the last term that
		\begin{align*}
			\osc_{\Omega_r}u&\int_{\Omega_r}|\nabla u||\nabla w|dx  
			\leq \zeta\int_{\Omega_r}|\nabla w|^2 dx + C_\zeta r^{4\beta}. 
		\end{align*}
		 Consequently, \eqref{eq:arkhi-previa-epsilon-zeta} leads us to the next inequality
	\[	(\sigma - C_2 \osc_{\Omega_r}u ) \int_{\Omega_r} \!\!|\nabla u|^2\min\{|w|, \ell\} {\vp dx}\\
				\leq  C_1\osc_{\Omega_r} u \int_{\Omega_r} |w|dx + \zeta\int_{\Omega_r}|\nabla w|^2dx + C_\zeta r^{4\beta}\]
		for all $\zeta > 0$.  Let us choose $r_0<\delta/2<1$ such that $C_2\osc_{\Omega_{r_0}}u \leq \sigma/2$ . Then, arguing as in \eqref{eq:arkhi-q-bounds} and using $r < 1$, we conclude that
		\begin{equation}
			\label{eq:arkhi-part-2}
			\frac{\sigma}{2}\int_{\Omega_r}|\nabla u|^2\min\{|w|, \ell\}{\vp dx}\leq 2\zeta\int_{\Omega_r}|\nabla w|^2dx + C_\zeta r^{4\alpha},
		\end{equation}
		for all $\zeta > 0$ and $r<r_0$.
		
		Combining \eqref{eq:arkhi-part-1} and \eqref{eq:arkhi-part-2} with $\zeta = \frac{\sigma^2}{16\ C_2}$, we have
		\begin{equation*}
			\frac{\sigma}{4}\int_{\Omega_r}|\nabla w|^2dx \leq  (C+\sigma)\int_{K_\ell}|\nabla w|^2dx + C r^{4\beta}
		\end{equation*}
		where $C$ does not depend on $\ell$. Thus, when $\ell\to \infty$ we obtain
		\begin{equation*}
			\int_{\Omega_r} |\nabla w|^2{\vp dx}\leq Cr^{4\beta} \quad \text{for }r \leq r_0 < \frac{\delta}{2}
		\end{equation*}
		Adapting \cite[Theorem 3.I]{Campanato1987} to our case, we deduce that
		\begin{equation*}
			\int_{\Omega_\rho}|\nabla v - (\nabla v)_{\overline x, \rho}|^2 dx \leq C\left(\frac{\rho}{r}\right)^{2+2\beta_0}\int_{\Omega_r}|\nabla v - (\nabla v)_{\overline x, r}|^2dx,
		\end{equation*}
		for some $\beta_0>0$ and $\rho < r$. Thus,
		\begin{align*}
			\int_{\Omega_\rho}|\nabla u - (\nabla u)_{\overline x, \rho}|^2dx &\leq 2\int_{\Omega_r}|\nabla w-(\nabla w)_{\overline x,r}|^2+2\int_{\Omega_r}|\nabla v-(\nabla v)_{\overline x,r}|^2\\
			& \leq C\rho^{4\beta}+C\left(\frac{\rho}{r}\right)^{2+2\beta_0}\int_{\Omega_r}|\nabla v - (\nabla v)_{\overline x, r}|^2dx \\ &\leq  C\rho^{4\beta}+ 2 C\left(\frac{\rho}{r}\right)^{2+2\beta_0}\left(\int_{\Omega_r}|\nabla u - (\nabla u)_{\overline x, r}|^2dx + r^{4\beta}\right).
		\end{align*}
		Set $\beta = \frac{1+\beta_0}{2}$. Accordingly, we get that
		\begin{equation*}
			\frac{1}{\rho^{2 +2\beta_0}}\int_{\Omega_\rho}|\nabla u - (\nabla u)_{\overline x, \rho}|^2 \, d x \leq 3C + \frac{2C}{r^{2 + 2\beta_0}} \int_{\Omega_r}|\nabla u - (\nabla u)_{\overline x, r}|^2dx,
		\end{equation*}
		so $\nabla u\in \mathcal L^{2,2 + 2\beta_0}(\Omega_{r/2})$. By \cite[Theorem I.2]{Campanato1963}, $\nabla u\in C^{0,\beta_0}(\Omega_{r/2})$. Therefore, we can conclude that $u \in C^{1, \beta_0}(\overline\Sigma)$. 
	\end{proof}
	
	\subsection{Higher regularity via continuity and bootstrapping methods} \label{boost}

In order to obtain a higher order of regularity, we need to smooth out  the function $f$ and to use classical results in regularity theory.
	
Since $f$ is a Lipchitz function, we can approximate it by a mollified sequence, which we call $f_\delta$, as in Section \ref{sec:notations}. Thanks to the extra regularity of this approximation, we can obtain the next statement.
	\begin{prop}\label{prop:holder-D2u}
		Let $u$ be a weak solution of \eqref{eq:relaxed-pde}. Then $u\in C^\infty(\overline\Sigma;\mathcal N)$.
	\end{prop}
\begin{proof}
To avoid technicalities, we assume that $\Sigma$ is an open bounded set in $\R^2$ (otherwise, we argue by partition of unity subordinate to a cover by isothermal coordinates, as usual).

We note that $u$ is a $C^{1,\beta_0}(\overline\Sigma; \R^N)$ solution to
\begin{equation}
  \label{schauder-1}
  \left\{\begin{array}
    {cc} \displaystyle{\rm div}\left(\left(\frac{\varrho}{\sqrt{\varepsilon^2\varrho^2+|\nabla w|^2}}+\sigma\right)\nabla w\right)=\Psi(x) & {\rm in \ }\Sigma \smallskip \\ \nabla w\cdot \nu=0 & {\rm on \ } \partial\Sigma,
  \end{array}\right.
\end{equation}
with $\displaystyle\Psi^\alpha(x):=-(b(\cdot,\nabla u)+\sigma)\sum_{i}\mathcal{A}_{u}^\alpha\left(\partial_i u, \partial_i u\right)-\lambda\varrho^2(\exp^{-1}_{u} f_\delta)^\alpha\in C^{0,\beta_0}(\overline \Sigma).$
Moreover, by standard comparison technique, we can show that the solution to \eqref{schauder-1} is unique.

We rewrite \eqref{schauder-1} as
\begin{equation}\label{schauder-2}
  \left\{\begin{array}
    {cc} \displaystyle Lw:=\sum_{i,j} A_{ij}(x,\nabla w)\partial^2_{ij}w+\Xi(x,\nabla w)=\Psi & {\rm in \ } \Sigma \\ \nabla w\cdot \nu=0 & {\rm on \ } \partial\Sigma,
  \end{array}\right.
\end{equation}
with $$A_{ij}(x,P)=\frac{\varrho}{\sqrt{\varepsilon^{2}\varrho^2 + |P|^{2}}} \left[ \delta_{ij}I^{N \times N} - \left( \frac{P_i}{\sqrt{\varepsilon^{2}\varrho^2 + |P|^{2}}}\otimes \frac{P_j}{\sqrt{\varepsilon^{2}\varrho^2 + |P|^{2}}}
 \right) \right] + \sigma \delta_{ij}I^{N \times N},$$

 $$\Xi[x,P]=\frac{\langle\nabla\varrho,P\rangle|P|^2}{(\varepsilon^2\varrho^2+|P|^2)^{\frac{3}{2}}}.$$

 It is easy to see that the system \eqref{schauder-2} satisfies the classical Agmon-Douglis-Nirenberg conditions \cite{agmon-douglis-nirenberg}. Note also that the system has at most one solution in $C^{2,\beta_0}(\overline\Sigma)$ due to uniqueness for \eqref{schauder-1}. Therefore, the following Schauder estimates hold for solutions  \cite[Remark 2]{agmon-douglis-nirenberg}:
 $$\|w\|_{C^{2,\beta_0}(\overline\Sigma)}\leq C \|\Psi\|_{C^{0,\beta_0}(\overline\Sigma)}. $$

 As before, one can check that $w=0$ is the unique solution to the  homogeneous Neumann problem for the system $$t\Delta w+(1-t) Lw=0.$$
Hence, by the continuity method (e.g. \cite[Theorem 5.25]{giacquinta-martinazzi}), we can prove that there exists a solution to \eqref{schauder-2} ${\rm w}\in C^{2,\beta_0}(\overline\Sigma)$. By uniqueness, $u={\rm w}$. Finally, a bootstrap technique shows that $u\in C^\infty(\overline\Sigma)$.
\end{proof}

	\section{Lipschitz regularity} \label{Lip}
	
		In this section and afterwards, we will use repeatedly the following notations:
	\[
	\text{s}_\kappa(t) = \begin{cases}
	\frac{\sin(\sqrt{\kappa}t)}{\sqrt{\kappa}}, &\kappa > 0,\\
	t. &\kappa = 0,\\
	\frac{\sinh(\sqrt{|\kappa|}t)}{\sqrt{|\kappa|}}, &\kappa < 0,
	\end{cases}  \qquad {\rm c}_\kappa = \text{s}_\kappa', \qquad {\rm ta}_\kappa = \frac{\text{s}_\kappa}{\text{c}_\kappa} \qquad \text{and} \quad \quad {\rm co}_\kappa = \frac{\text{c}_\kappa}{\text{s}_\kappa}.
	\]
	
	\subsection{Energy density as a subsolution of an elliptic PDE}
	
	As in \cite{Porr} we get Lipschitz regularity for general (non-necessarily convex) domains, but we have the extra difficulty of dealing with maps into a target manifold (instead of scalar-valued) and on a curved domain. Typically, regularity follows from integral estimates as in \cite{GiacomelliLasicaMoll2019}, but with the drawback of having to require convexity of the boundary. To overcome this, we adapt to our setting the Bernstein technique used in \cite{Porr}; let us highlight that this approach only works arguing with the intrinsic version of our elliptic system  {\normalfont($\mathcal S_{\varepsilon, \sigma}^{f_\delta}$)}, namely, set $w:= |d u|^2$ where $u:= u_{\varepsilon, \sigma}^\delta$ solves the system
	\begin{equation}\label{relaxed-pde-intr}
	\left\{ \begin{array}{rcll}
\ds	\dv_g \bigg[\Big(\frac1{\sqrt{w + \varepsilon^2}} + \sigma \Big) du\bigg] & = & -\lambda\exp^{-1}_u f_\delta &\qquad \text{in $\Sigma$} \medskip \\
	\nu\cdot  d u& = & 0 &\qquad \text{on $\partial\Sigma$}
	\end{array}\right.\tag{$\mathcal S_{\varepsilon, \sigma}^{f_\delta}$}
	\end{equation}
with $f_\delta$ the mollified version of the initial $f$ defined as in \eqref{moli-f}.

The goal is to get uniform bounds for $w$, which are independent of the parameters $\varepsilon, \sigma$ and $\delta$ in order to allow a limiting process when any of these values goes to 0. As a prior step, we derive an elliptic PDE for $w$ to which we will apply a maximum principle argument.
With this aim, first notice that the system at the first line in \eqref{relaxed-pde-intr} can be rewritten as
\begin{equation} \label{eqn-intr-2}
\bigg(\frac1{\sqrt{w + \varepsilon^2}} + \sigma \bigg) \tau(u)^\alpha \,  - \frac1{2(w +\varepsilon^2)^{3/2}} \<\nabla w, \nabla u^\alpha\>_g = \, -\lambda\big(\exp^{-1}_u f_\delta\big)^\alpha
\end{equation}
for $\alpha = 1, \ldots, n$	and $\tau$ be the 2-tension (recall Section \ref{tensions}).

 We define the operator
\begin{equation} \label{def-L}
\mathscr L w:= - \bigg(\frac1{\sqrt{w + \varepsilon^2}}  + \sigma \bigg) \Delta_g w + \frac1{(w +\varepsilon^2)^{3/2}} (h_{\alpha \beta} \circ u) \nabla^2 w(\nabla u^\alpha, \nabla u^\beta).
\end{equation}
Observe that taking normal coordinates around a fixed but arbitrary point $x \in \Sigma$ and corres\-pondingly around $u(x)\in \mathcal N$, it is easy to check $\mathscr L$ is an elliptic operator.
\begin{lemma} \label{Porr:lem1}
 Suppose that the supremum of the sectional curvatures of $\mathcal N$ within $B_h(p,R)$ is non-positive. Then $w$ as above satisfies
\begin{align*}
     \mathscr L w + \lambda w  + 2 \bigg(\frac1{\sqrt{w + \varepsilon^2}}  + \sigma \bigg)\bigg(\big|\nabla^{g\boxtimes h} d u\big|^2_{g\boxtimes h} - \mathscr R  \bigg)  \leq & \frac3{2} \frac{\big<dw \otimes \<\nabla w, \nabla u\>_g, du\big>_{g\boxtimes h}}{(w +\varepsilon^2)^{5/2}} + \lambda  |df_\delta|^2 \\ & -  \frac{\<dw \otimes \tau(u), du\>_{g\boxtimes h}}{(w +\varepsilon^2)^{3/2}}   - \frac{|\nabla w|^2_g}{2(w +\varepsilon^2)^{3/2}},
\end{align*}
where the curvature term $\mathscr R$ is given by
\begin{equation} \label{curv-BW}
\mathscr R:= \sum_{i.j =1}^2 \big<\mathcal R^{u^\ast h}(e_i, e_j) e_j, e_i\big>_{u^\ast h} - {\rm Ric}_\Sigma(du, du), \quad \text{with} \quad \mathcal R^{u^\ast h}(X, Y) Z := \mathcal R^h(u_\ast X, u_\ast Y) u_\ast Z,
\end{equation}
for a local orthonormal frame $\{e_i\}_{i =1,2}$ of $\Sigma$.
\end{lemma}	
\begin{proof}
Recall that on $T^\ast \Sigma \otimes u^\ast T \mathcal N$ we have the natural bundle metric $g\boxtimes h:= g^{-1} \otimes u^\ast h$ so that $\<A, B\>_{g\boxtimes h}:= g^{ij} (h_{\alpha \beta} \circ u) A_i^\alpha B_j^\beta$, while $\widetilde \nabla$ denotes the linear connection on $u^\ast T\mathcal N$ given by \eqref{ind-conn}. Differentiating \eqref{eqn-intr-2} and multiplying with $du$, using the inner product induced by $g\boxtimes h$, we reach
\begin{align*}
\bigg(\frac1{\sqrt{w + \varepsilon^2}} &+ \sigma \bigg)   \big<\widetilde \nabla \tau(u), du\big>_{g\boxtimes h} \,  - \frac1{2} \frac{\<dw \otimes \tau(u), du\>_{g\boxtimes h}}{(w +\varepsilon^2)^{3/2}} + \frac3{4} \frac{\big<dw \otimes \<\nabla w, \nabla u\>_g, du\big>_{g\boxtimes h}}{(w +\varepsilon^2)^{5/2}} \\ & - \frac1{2 (w +\varepsilon^2)^{3/2}} \Big((h_{\alpha \beta} \circ u) \nabla^2 w(\nabla u^\alpha, \nabla u^\beta) + \frac1{2}|\nabla w|^2_g\Big)= \, -\lambda\big<\widetilde \nabla  \big(\exp^{-1}_u f_\delta\big), du\big>_{g\boxtimes h}.
\end{align*}
Replacing the first term on the left hand side by means of Bochner-Weitzenb\"ock formula (see e.g. \cite[p. 128]{Nishi}), we have
\begin{align} \label{Porr-1}
\bigg(\!\frac1{\sqrt{w + \varepsilon^2}}  &+ \sigma \!\bigg) \bigg(\frac{\Delta_g w}{2} - \big|\nabla^{g\boxtimes h} d u\big|^2_{g\boxtimes h} + \mathscr R \! \bigg)    -  \frac{\<dw \otimes \tau(u), du\>_{g\boxtimes h}}{2 (w +\varepsilon^2)^{3/2}} + \frac3{4} \frac{\big<dw \otimes \<\nabla w, \nabla u\>_g, du\big>_{g\boxtimes h}}{(w +\varepsilon^2)^{5/2}} \nonumber \\ & \hspace*{-0.3cm} - \frac1{2 (w +\varepsilon^2)^{3/2}} \Big((h_{\alpha \beta} \circ u) \nabla^2 w(\nabla u^\alpha, \nabla u^\beta) + \frac1{2}|\nabla w|^2_g\Big)= \, -\lambda\big<\widetilde \nabla  \big(\exp^{-1}_u f_\delta\big), du\big>_{g\boxtimes h}
\end{align}
where $\mathscr R$ is defined as in \eqref{curv-BW}.

On the other hand, for any $X \in T\Sigma$, by the chain rule we have
\begin{align*}
\widetilde \nabla_X  \big(\exp^{-1}_u f_\delta\big) = ^h\!\nabla_{u_\ast X}  \big(\exp^{-1}_\cdot f_\delta\big) + D_{{f_\delta}_\ast X} \exp^{-1}_u \cdot.
\end{align*}
Here notice that for each $x \in \Sigma$, $D$ stands for the standard directional derivative of the vector space $T_{u(x)} \mathcal N$. Hence, for the term on the right hand side of \eqref{Porr-1}, we can write
\begin{align*}
- \big<\widetilde \nabla  \big(\exp^{-1}_u f_\delta\big), du\big>_{g\boxtimes h} & = -\sum_{i =1}^2 \Big( \big<^h\!\nabla_{u_\ast e_i}  \big(\exp^{-1}_\cdot f_\delta\big), u_\ast e_i\big>_h  +\big<D_{{f_\delta}_\ast e_i} \exp^{-1}_u \cdot, u_\ast e_i\big>_h\Big)
\\ & = \nabla^2\bigg(\frac{d_h^2(\cdot,f_\delta)}{2}\bigg)(du, du) -\sum_{i =1}^2 \big<D_{{f_\delta}_\ast e_i} \exp^{-1}_u \cdot, u_\ast e_i\big>_h
\\ & \geq w  -\sum_{i =1}^2 \big<D_{{f_\delta}_\ast e_i} \exp^{-1}_u \cdot, u_\ast e_i\big>_h,
\end{align*}
which follows from the Hessian comparison theorem for the squared distance function for $\kappa \leq 0$ (see e.g. (10.47) in \cite[p. 266]{Villani}). Moreover, the first equality is obtained as in \cite[(7)]{CRMP}.

Now, to further estimate the second term in the right hand side, consider the geodesic on $\mathcal N$ given by  $\ga_i(t):= \exp_{f_\delta(x)} (t f_{\delta \ast x} e_i)$. Then it is well-known (see, for instance, \cite[p. 540]{Karcher1997}) that
\[\frac{d}{dt}\bigg|_{t=0} \exp_{u(x)}^{-1}\big(\gamma_i(t)\big) = J'(0) \in T_{u(x)}\mathcal N,\]
where $J(s)$ is the Jacobi field along the geodesic ${\vp \zeta}(s) = \exp_{u(x)}(s \exp_{u(x)}^{-1} f_\delta(x))$ determined by $J(0) = 0$ and $J(1) = \ga_i'(0) = f_{\delta \ast x} e_i \in T_{f_\delta(x)} \mathcal N$ (see Figure \ref{fig:Jac-dibu}). Then one can reproduce almost verbatim the proof of Jacobi field comparison in \cite[Theorem 11.2]{Lee-old}, but for a non-unit speed geodesic, as in our case $|{\vp \zeta}'|_h = d_h(u(x), f_\delta(x))$, and a non-necessarily normal Jacobi field, to deduce that $t |J'(0)| \leq |J(t)|$. In particular, for $t = 1$, we reach
\[|J'(0)| \leq |J(1)| =\big|f_{\delta \ast x} e_i\big|_h. \]

\begin{figure}[h]
	\includegraphics[scale =0.6]{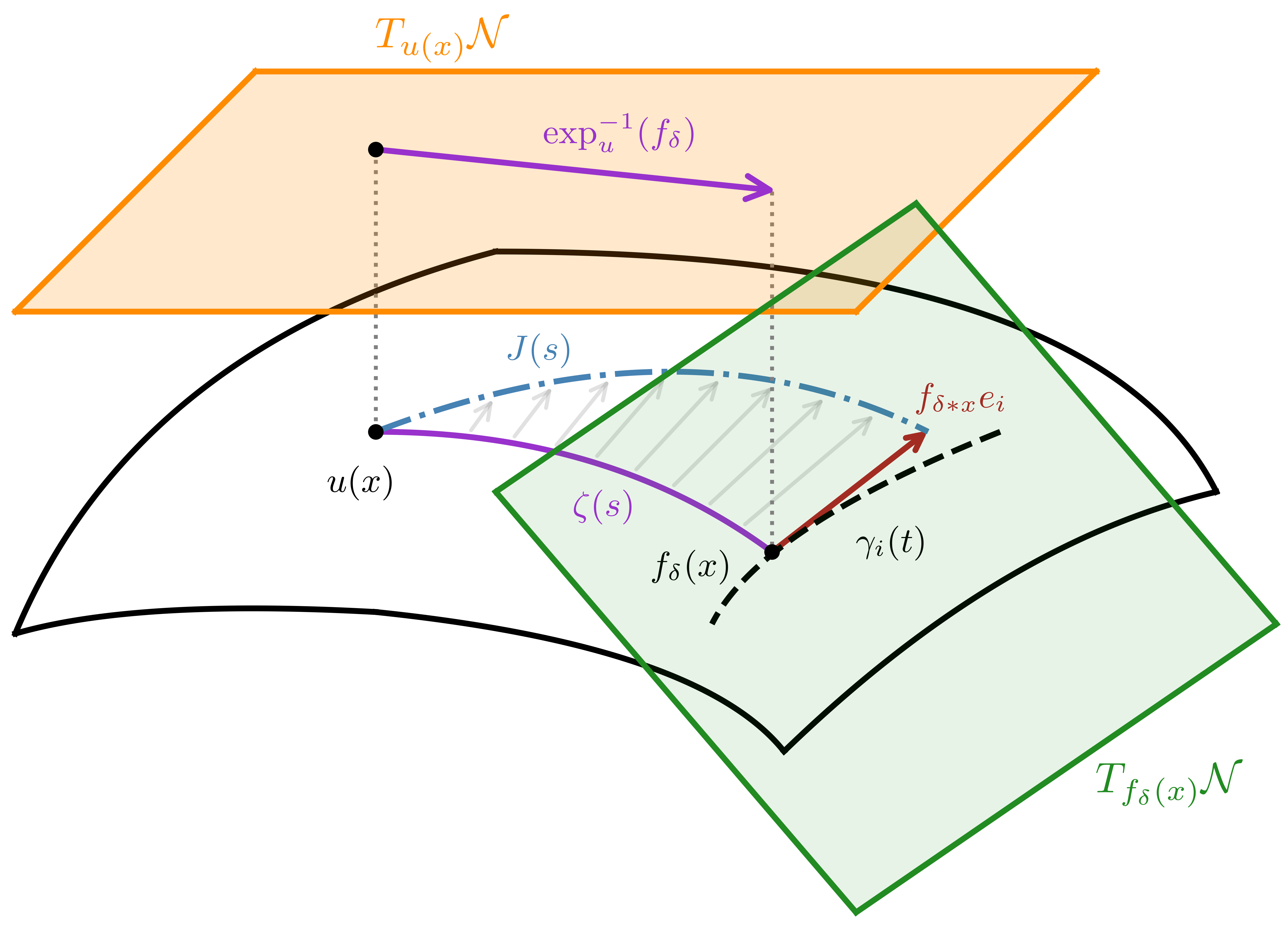}	
	\vspace*{-0.3cm}\caption{Setting to apply Jacobi field comparison.} \label{fig:Jac-dibu}	
\end{figure}

Taking all the above considerations into account and using Young's inequality, we get
\begin{align*}
- \big<\widetilde \nabla  \big(\exp^{-1}_u f_\delta\big), du\big>_{g\boxtimes h} &  \geq w  - |df_\delta| \cdot |du| \geq \frac1{2} \big(w - |df_\delta|^2\big).
\end{align*}

Substituting into \eqref{Porr-1} and rearranging terms, we can write
\begin{align*}
      \frac3{4} \frac{\big<dw \otimes \<\nabla w, \nabla u\>_g, du\big>_{g\boxtimes h}}{(w +\varepsilon^2)^{5/2}} -  \frac{\<dw \otimes \tau(u), du\>_{g\boxtimes h}}{2 (w +\varepsilon^2)^{3/2}}   &- \frac{|\nabla w|^2_g}{4 (w +\varepsilon^2)^{3/2}}    \geq  \, \frac{\lambda}{2} \big(w - |df_\delta|^2\big) \\ &  \hspace*{-0.6cm} + \frac1{2}\mathscr L w + \bigg(\frac1{\sqrt{w + \varepsilon^2}}  + \sigma \bigg)\bigg(\big|\nabla^{g\boxtimes h} d u\big|^2_{g\boxtimes h} - \mathscr R  \bigg).
\end{align*}
Further rearrangement of the terms and multiplication by 2 leads to the inequality in the statement.
\end{proof}

Let us remark that, as $\Sigma$ is a compact surface, there is always a constant $C_{_\Sigma}$ so that ${\rm Ric}_\Sigma \geq - C_{_\Sigma} g$,  which hence does not introduce any further restriction in our current setting. We keep the subindex to identify this constant. We conclude that $w$ is a subsolution of an elliptic equation. More precisely,
\begin{cor} \label{Porr-Cor1}	Assume that the upper bound $\kappa$ for the sectional curvatures of $\mathcal N$ within $B_h(p, R)$ is non-positive. Then $w = |d u|^2$, for $u$ the unique solution of \eqref{relaxed-pde-intr}, satisfies the inequality
\begin{align*}
\mathscr L w + \lambda w  -2  \bigg(\frac1{\sqrt{w + \varepsilon^2}}  + \sigma \bigg) C_{_\Sigma} w 
\leq \frac{3}{2} \frac{|\nabla w|^2_g}{(w +\varepsilon^2)^{3/2}} + \lambda  |df_\delta|^2.
\end{align*}	
\end{cor}
\begin{proof}
Notice that, as $\nabla^{g\boxtimes h} d u$	is a symmetric $T \mathcal N$-valued 2-tensor, for any orthonormal basis $\{e_1, e_2\}$ of $\Sigma$, we have
\[\big|\nabla^{g\boxtimes h} d u\big|^2_{g\boxtimes h} = \sum_{i,j=1}^2 \big|\nabla^{g\boxtimes h} d u(e_i, e_j)\big|^2_h \geq \frac1{2} \sum_{i=1}^2 \Big|\nabla^{g\boxtimes h} d u(e_i, e_i)\Big|^2_h = \frac1{2} |\tau(u)|^2_h. \]
On the other hand,  by means of Cauchy-Schwarz and Young's inequality, we can estimate
\begin{align*}
 \frac{\<dw \otimes \tau(u), du\>_{g\boxtimes h}}{(w +\varepsilon^2)^{3/2}} & \leq  \frac{|dw \otimes \tau(u)|_{g\boxtimes h} \, \sqrt{w}}{(w +\varepsilon^2)^{3/2}} \leq \frac{|\tau(u)|_{h}}{(w +\varepsilon^2)^{1/4}}\cdot \frac{|dw|_g   \, \sqrt{w}}{(w +\varepsilon^2)^{5/4}} \\ & \leq \frac{|\tau(u)|_{h}^2}{2(w +\varepsilon^2)^{1/2}}   + \frac{|\nabla w|_g^2   \, w}{2 (w +\varepsilon^2)^{5/2}} \leq \frac{\big|\nabla^{g\boxtimes h} d u\big|^2_{g\boxtimes h}}{(w +\varepsilon^2)^{1/2}}   + \frac{|\nabla w|_g^2 }{2 (w +\varepsilon^2)^{3/2}}
\end{align*}

Similarly, for the first term on the right hand side of the inequality in Lemma \ref{Porr:lem1}, we get
\begin{align} \label{term52}
 \frac{\big<dw \otimes \<\nabla w, \nabla u\>_g, du\big>_{g\boxtimes h}}{(w +\varepsilon^2)^{5/2}} \leq \frac{|\nabla w|_g^2 \, w }{ (w +\varepsilon^2)^{5/2}} \leq \frac{|\nabla w|_g^2 }{ (w +\varepsilon^2)^{3/2}}.
\end{align}

Accordingly, we reach
\begin{align*}
\mathscr L w + \lambda w  + 2 \bigg(\frac1{\sqrt{w + \varepsilon^2}}  + \sigma \bigg)\bigg(\big|\nabla^{g\boxtimes h} d u\big|^2_{g\boxtimes h} - \mathscr R  \bigg)  \leq &  \frac3{2} \frac{|\nabla w|^2_g}{(w +\varepsilon^2)^{3/2}} + \lambda  |df_\delta|^2 + \frac{\big|\nabla^{g\boxtimes h} d u\big|^2_{g\boxtimes h}}{(w +\varepsilon^2)^{1/2}}.
\end{align*}
Finally, discarding the non-negative term $\Big(\frac1{\sqrt{w +  \varepsilon^2}}  + 2 \sigma\Big) \big|\nabla^{g\boxtimes h} d u\big|^2_{g\boxtimes h}$ which results on the left hand side, and using the curvature assumptions, we deduce the inequality in the statement.
\end{proof}

	\subsection{Baby case: Lipschitz regularity in the convex case}

At this stage, let us point out that the above estimates can be applied directly to provide a  uniform bound on the Lipschitz constants of the approximating solutions
in the convex case. With this goal, let us denote by $r_{_\partial} :=d_g(\cdot, \partial \Sigma)$ the Riemannian distance to the boundary of our domain surface. As $\partial \Sigma$ is smooth, the function $r_{_\partial}$ is known to be smooth in the interior of $\Sigma$ out of the cut locus of $\partial \Sigma$ (see \cite[Proposition 3.10]{Sakurai}). In addition, $\nabla r_{_\partial}{|_{\partial\Sigma}} = -\nu$, where $\nu$ is the unit outward normal to $\partial \Sigma$, and hence for any $X \in T \partial \Sigma$ we have
\begin{equation} \label{dist-curv}
\nabla^2 r_{_\partial}(X, X) = \<\nabla_X \nabla r_{_\partial}, X\>_g = - \<\nabla_X \nu, X\>_g = - k_{_\partial} |X|_g^2,
\end{equation}
being $k_{_\partial}$ the geodesic curvature of the boundary curve $\partial \Sigma$. Consequently, for a convex boundary, the distance function $r_{_\partial}$ is concave.

\begin{prop} \label{Porr-cvx}
	Let $\mathcal N$ be a complete $n$-manifold with a non-positive upper bound $\kappa$ for the sectional curvatures within $B_h(p, R)$. Assume that $f$ is a Lipschitz function with $f(\overline\Sigma) \subset B_h(p, R)$ for some $R < R_\kappa$. Then  the unique solution  $u$ of \eqref{relaxed-pde-intr}
is  Lipschitz with constant only depending on the Lipschitz constant of $f$, under the extra assumption that the boundary curve $\partial \Sigma$ is convex.
\end{prop}
\begin{proof}
	A routine computation in local coordinates gives
	\[\partial_k w = 2 g^{ij} h_{\alpha \beta}|_u \Big(\partial_{ik}^2 u^\alpha - \partial_\ell u^\alpha \, \Gamma^{\ell}_{ik} + ^h\!\Gamma_{\delta \gamma}^\alpha\big|_u \,\partial_i u^\gamma \,\partial_k u^\delta\Big) \partial_j u^\beta,\]
where the shortcut $F|_u$ means composition with the map $u$.

Now choose a chart $\{x^1, x^2\}$ so that $x^2 = r_{_\partial}$, and thus for the normal derivative we have $\nabla_\nu = - \partial_2$. By normalization of the first vector, we can assume that $\{\partial_1, \partial_2\}$ is a local orthonormal frame. Then, taking $k =2$ in the above formula and using the Neumann boundary condition to remove the last term, we reach
\begin{align*}
\frac{\partial w}{\partial \nu} = -\partial_2 w &= - 2 g^{ij} h_{\alpha \beta}|_u  \nabla^2_{i2} u^\alpha \partial_j u^\beta = -2 h_{\alpha \beta}|_u \nabla^2_{12} u^\alpha \partial_1 u^\beta.
\end{align*}
On the other hand, notice that
\begin{align*}
\nabla_{12}^2 u^\alpha &= \<\nabla_{\partial_1} \nabla u^\alpha, \partial_2\> = \partial_1\<\nabla u^\alpha, \partial_2\> - \<\nabla u^\alpha, \nabla_{\partial_1} \partial_2\>.
\end{align*}
Moreover, by the Neumann boundary condition, $\<\nabla u^\alpha, \partial_2\>$ vanishes on the boundary, and thus its gradient points in the normal direction. In addition, $\nabla u^\beta = \partial_1 u^\beta \, \partial_1 \in T \partial \Sigma$. Accordingly, we can write
\begin{align} \label{Hopf}
\frac{\partial w}{\partial \nu} = -\partial_2 w & =-2 h_{\alpha \beta}|_u \<\nabla \<\nabla u^\alpha, \partial_2\>, \nabla u^\beta\> + 2 h_{\alpha \beta}|_u \partial_1 u^\alpha \partial_1 u^\beta \<\partial_1, \nabla_{\partial_1} \partial_2\> \nonumber
\\ & = 0 - w \<\partial_1, \nabla_{\partial_1} \nu\> = - k_{_\partial} w \leq 0,
\end{align}
which follows by applying the convexity hypothesis.

The latter estimate, by means of Hopf's maximum principle, ensures that the maximum of $w$ cannot be attained at the boundary $\partial \Sigma$. Let us then evaluate the inequality from Corollary \ref{Porr-Cor1} at a maximum interior point $x$, taking into account that $\nabla w = 0$ and $\mathscr L w \geq 0$. Then discarding also the squared norm on the left hand side, we reach
\begin{align*}
 (\lambda   -2\sigma C_{_\Sigma}) w \leq 2 \frac1{\sqrt{w + \varepsilon^2}}  C_{_\Sigma} w  + \lambda  |df_\delta|^2  \leq 2   C_{_\Sigma} \sqrt{w}  + \lambda  |df_\delta|^2.
\end{align*}	
Without loss of generality, we can take $\sigma$ small enough so that $\sigma C_{_\Sigma} < \lambda/4$ and assume that $w > 1$ (otherwise we are done). This leads to
\[ \frac{\lambda}{2} \sqrt{w} \leq 2   C_{_\Sigma} + \frac{\lambda}{\sqrt{w}}  |df_\delta|^2 \leq 2   C_{_\Sigma} + \lambda C \, {\rm Lip}(f)^2,\]
which implies the stated Lipschitz regularity of $u$ in $\overline \Sigma$.
\end{proof}

	\subsection{Modification of the argument for a weighted energy density} \label{peso}
To achieve our main result removing the extra condition about convexity of the boundary, one needs to work harder. The key idea comes from \cite[Lemma 2.4]{LeoPorr} and amounts to multiply $w$ by a suitable test function so that for the product the nonnegativity in \eqref{Hopf} still holds and we can argue at a maximum interior point with an appropriate version of Corollary \ref{Porr-Cor1}, which we deduce in the following lines.
\begin{lemma} \label{Porr-Lem2}
Set $\widetilde w:= \varphi \, w$ for some positive function $\varphi \in C^2(\overline \Sigma)$. Then, under the same assumptions and notation as in Lemma \ref{Porr:lem1} and at a maximum interior point of $\Sigma$, it holds
\begin{align*}
\lambda \widetilde w  -2  \bigg(\frac1{\sqrt{w + \varepsilon^2}}  + \sigma \bigg) & C_{_\Sigma} \widetilde w   \leq    \lambda \, \varphi  |df_\delta|^2 + \sigma w \bigg(2 \frac{|\nabla \varphi|_g^2}{\varphi} -\Delta \varphi\bigg) + \frac1{2} \sqrt{ w} \,\bigg(11 \frac{|\nabla \varphi|_g^2}{\varphi} + 5 |\nabla^2 \varphi|_g   \bigg).
\end{align*}
\end{lemma}
\begin{proof}
Taking into account that the Hessian of a product reads as
\[\nabla^2_{ij} \widetilde w = \varphi \,\nabla^2_{ij} w + \partial_i w \, \partial_j \varphi + \partial_j w \, \partial_i \varphi + w \, \nabla^2_{ij}\varphi,\]
we can write
\begin{align}
\mathscr L \widetilde w  = & \varphi \, \mathscr Lw -  \bigg(\frac1{\sqrt{w + \varepsilon^2}}  + \sigma \bigg) \Big(w \, \Delta \varphi + 2 \<\nabla w, \nabla \varphi\>_g\Big) + \frac{w}{(w+ \varepsilon^2)^{3/2}} (h_{\alpha \beta}\circ u) \nabla^2 \varphi(\nabla u^\alpha, \nabla u^\beta) \\ & + \frac{2}{(w+ \varepsilon^2)^{3/2}} (h_{\alpha \beta}\circ u) \<\nabla w, \nabla u^\alpha\>_g \big<\nabla \varphi, \nabla u^\beta\big>_g  \label{Ltest}
\\ & \leq \varphi \, \mathscr Lw -  \bigg(\frac1{\sqrt{w + \varepsilon^2}}  + \sigma \bigg) \Big(w \, \Delta \varphi + 2 \<\nabla w, \nabla \varphi\>_g\Big) + \frac{|\nabla^2 \varphi|_g w}{(w+ \varepsilon^2)^{3/2}} + 2  \frac{|\nabla w|_g |\nabla \varphi|_g w}{(w+ \varepsilon^2)^{3/2}}, \nonumber
\end{align}
which follows by applying Cauchy-Schwarz arguing with normal coordinates both in the domain and target, around a point $x \in \Sigma$ and $u(x) \in \mathcal N$, respectively.

Now multiply by $\varphi$ the inequality in Corollary \ref{Porr-Cor1} and substitute the above estimate to reach
\begin{align*}
\mathscr L \widetilde w + \lambda \widetilde w  -2  \bigg(\frac1{\sqrt{w + \varepsilon^2}}  + \sigma \bigg) & C_{_\Sigma} \widetilde w   \leq   \,  \frac{3 \varphi |\nabla w|^2_g}{2 (w +\varepsilon^2)^{3/2}} + \lambda \, \varphi  |df_\delta|^2 -  w \bigg(\frac1{\sqrt{w + \varepsilon^2}}  + \sigma \bigg)  \Delta \varphi \\ & + 2 \bigg(\frac1{\sqrt{w + \varepsilon^2}}  + \sigma \bigg) |\nabla w|_g |\nabla \varphi|_g + |\nabla^2 \varphi|_g  \sqrt{w} + 2  \frac{|\nabla w|_g |\nabla \varphi|_g}{\sqrt{w+ \varepsilon^2}}.
\end{align*}

Now if we compute the previous inequality at a maximum interior point $x \in \Sigma$, notice that
\begin{equation} \label{grad-max}
0 = \nabla \widetilde w = (\nabla \varphi) w + \varphi \nabla w, \quad \text{and hence } \quad \nabla w = - \frac{\nabla \varphi}{\varphi} w.
\end{equation}
 By substitution of the latter, we get
\begin{align*}
\mathscr L \widetilde w + \lambda \widetilde w  -2  \bigg(\frac1{\sqrt{w + \varepsilon^2}}  + \sigma \bigg) & C_{_\Sigma} \widetilde w   \leq   \,  \frac{3 |\nabla\varphi|^2_g}{2 \varphi} \sqrt{w}+ \lambda \, \varphi  |df_\delta|^2 - \sigma w \Delta \varphi + \sqrt{2 \, w} \, |\nabla^2 \varphi|   \\ & + 2 \bigg(\frac1{\sqrt{w + \varepsilon^2}}  + \sigma \bigg)  \frac{|\nabla \varphi|_g^2}{\varphi} w + |\nabla^2 \varphi|_g  \sqrt{w} + 2  \frac{|\nabla\varphi|^2_g}{\varphi} \sqrt{w},
\end{align*}
where we have also used that $(\Delta \varphi)^2 \leq 2 |\nabla^2 \varphi|^2$. Rearranging the right hand side and taking into account that $\mathscr L \widetilde w \geq 0$ at a maximum interior point, we deduce the estimate in the statement.
\end{proof}

	\subsection{Lipschitz regularity of the approximated minimizers}
		\begin{prop}\label{prop:lipschitz}
			Let $u$ be the minimizer of $\mathcal E_{\sigma, \varepsilon}$ when $\kappa \leq 0$. Then there exist constants $C >0$ and  $\widetilde{C}_{_\Sigma}$, depending only on the domain surface $\Sigma$, such that the following uniform bound
			\[
			|d u| \leq C {\rm Lip} (f)^2 + \widetilde{C}_{_\Sigma}
			\]
		holds. 
	\end{prop}
	\begin{proof}
	By the already mentioned regularity properties of the distance to the boundary $r_{_\partial}$, we can find a constant $\ell > 0$ so that $r_{_\partial}$ is $C^\infty$ on the set $\{x \in \Sigma \, : \, r_{_\partial}(x) \leq \ell\}$. Then consider a function $\eta$ which coincides with $r_{_\partial}$ if $r_{_\partial} \leq \ell/2$, and vanishes for $r_{_\partial} \geq \ell$.
	
		Let us take in particular $\widetilde w := e^{c \, \eta} w$ with $c := 2 \max_{\partial \Sigma} |k_{_\partial}|$. With the same conventions as in the proof of Proposition \ref{Porr-cvx} and using \eqref{Hopf}, on boundary points we have
	\begin{align} \label{bdry-ineq2}
	\frac{\partial \widetilde w}{\partial \nu} = e^{c \, \eta} w \big(- k_{_\partial} - c \big) \leq \widetilde w (|k_{\partial}| - c) \leq 0,
	\end{align}
	being the latter true by the choice of $c$. Consequently, arguing as in the proof of Proposition \ref{Porr-cvx}, Hopf's maximum principle ensures that the maximum of $\widetilde w$ occurs at the interior of $\Sigma$. By composing $r_{_\partial}$ with any cut-off function in $W^{2, \infty}$, we have the following growth for the derivatives of our concrete test function $\varphi = e^{c \, \eta}$:
	\begin{equation} \label{bound-test}
	\frac{|\nabla \varphi|^2}{\varphi} + \frac1{2} |\nabla^2 \varphi| \leq \frac{3}{2} c^2 \varphi |\nabla \eta|^2 + \frac{c}{2} \varphi |\nabla^2 \eta| \leq C \varphi \bigg(\frac{3}{2} c^2 + \frac{c^2}{4}\bigg) \leq \widetilde{C}_{_\Sigma} \varphi,
	\end{equation}
	where $C$ comes from the bounds of the first and second derivatives of the cut-off function and we have also used \eqref{dist-curv}. The notation indicates that the constant only depends on $\Sigma$, but is uniform in $\sigma$ and $\varepsilon$, and its concrete meaning may change from line to line.
	
	We are now in position to apply Lemma \ref{Porr-Lem2} to reach
	\begin{align*}
	(\lambda -  2\sigma  C_{_\Sigma}) \widetilde w  \leq 2  \, C_{_\Sigma} \varphi {\sqrt w}   +  \lambda \, \varphi  |df_\delta|^2 + (\sigma w + \sqrt{w}) \widetilde{C}_{_\Sigma} \varphi,
	\end{align*}
	and rearranging terms leads to
	\[\big(\lambda -  2\sigma  C_{_\Sigma} -\sigma \widetilde{C}_{_\Sigma}\big) \sqrt{w}   \leq 2  \, C_{_\Sigma}    +  \frac{\lambda}{\sqrt{w}} \,  |df_\delta|^2 +  \widetilde{C}_{_\Sigma},\]
		Finally, taking $\sigma$ small enough so that $\sigma(2  C_{_\Sigma} -\sigma \widetilde{C}_{_\Sigma}) < \frac{\lambda}{2}$ and arguing as at the end of the proof of Proposition \ref{Porr-cvx}, we deduce the estimate and hence the regularity claim in the statement. In fact, notice that we have a bound for $\max \widetilde w$, but the conclusion follows because $|d u|^2 \leq \widetilde w$.
	\end{proof}	

	\subsection{Convergence and Lipschitz regularity of the original problem}\label{sect:conv}

	\begin{prop}\label{final_result}
		Let $\sigma \geq 0$ and $u$ be the minimizer of $\mathcal E_\sigma$ when $\kappa \leq 0$. Then $u\in C^{0,1}(\overline\Sigma;\mathcal N)$.
	\end{prop}
	\begin{proof}

	Let $\{\varepsilon_k\}_k,\{\delta_k\}_k$ be such that $\varepsilon_k,\delta_k\to 0$, and let $u_k$ be a weak solution of {\normalfont($\mathcal S_{\varepsilon_k, \sigma}^{f_{\delta_k}}$)}. We can keep $\sigma>0$ fixed, because for the case when $\sigma=0$, we simply take $\sigma_k\to 0$ and repeat all the computations.  By Propositions \ref{prop:holder-D2u} and \ref{prop:lipschitz}, we have that the functions $u_k$ are uniformly bounded in $C^{0,1}(\overline\Sigma;\mathcal N)$.

Therefore, one can find a subsequence (which will not be relabeled for brevity) such that $u_k\to u$ strongly in $L^\infty(\overline\Sigma;\mathcal N)$ to some $u\in C^{0,1}(\overline\Sigma;\mathcal N)$. Moreover, working again in isothermal coordinates for the domain and in extrinsic by isometric embedding in $\R^N$, we have that $\nabla u_k{\overset{*}{\rightharpoonup}} \nabla u$ in $L^\infty(\Omega_\ell)$ for all $\ell$.

	Remember from \eqref{eq:simplified-pde}, that it holds $$a_{\varepsilon_k}=\left(\frac{\varrho}{\big(|\nabla u_k|^2+\varepsilon_k^2\varrho^2\big)^{1/2}}+\sigma\right)\nabla u_k=:\varrho \, B_k+\sigma\nabla u_k.$$

Since $\|B_k||_\infty\leq 1$,  $$B_k{\overset{*}{\rightharpoonup}}\, ^\ell B\ {\rm\ in \ }L^\infty(\Omega_\ell) {\rm \ and \ }\| ^\ell B\|_\infty\leq 1,$$ for a subsequence, not relabeled.
Thus	

\[
	a_{\varepsilon_k} {\overset{*}{\rightharpoonup}}\,  ^\ell Z=\varrho\,\, ^\ell B+\sigma \nabla u \quad \text{in $L^\infty(\Omega_\ell)$}.
	\]

	From \eqref{eq:simplified-pde}, we deduce that $\dv(a_{\varepsilon_k})$ is uniformly bounded in $L^2(\Omega_\ell)$. Applying \cite[\S 5, Lemma 4]{Evans1988}, we deduce that
	\[
	\langle \nabla u_k,   a_{\varepsilon_k}\rangle \rightharpoonup \langle \nabla u , \varrho\, \, ^\ell B\rangle+\sigma|\nabla u|^2 \quad \text{in $L^2(\Omega_\ell)$.}
	\]
	Furthermore, from the definition of $B_k$ we note that $\langle \nabla u_k, B_k\rangle \geq |\nabla u_k| - \varepsilon_k$ and by the lower semicontinuity of the modulus with respect to weak convergence, we have that
	\begin{equation}
		\label{eq:Du_inequality}
		\varrho|\nabla u|+\sigma|\nabla u|^2\leq\liminf_{k\to\infty} \langle \nabla u_k,   a_{\varepsilon_k}\rangle\leq \varrho\langle \, ^\ell B,\nabla u\rangle+\sigma|\nabla u|^2,
	\end{equation}
which implies that $\langle ^\ell B,\nabla u\rangle=|\nabla u|$.

Moreover, by the strong convergence to $u$ we obtain that $ ^\ell Z(x)\in T_{u(x)}\mathcal N$, for a.e. $x\in\Omega_\ell$, which in turn leads to $ ^\ell B(x)\in T_{u(x)}\mathcal N$, for a.e. $x\in\Omega_\ell$.

We let now $\mathcal B:=\sum_\ell  \chi_\ell \cdot (^\ell B \circ \phi_\ell)$ and define the 1-form $\mathcal Z$ as $$\mathcal Z (\partial_i):=\mathcal B_i.$$ Then,  we deduce that $\mathcal Z$ satisfies
\eqref{eq:def-Z}.

On the other hand, by  strong convergence and  boundedness of $\dv(\mathcal Z_{\varepsilon_k})$, we can take limits in the weak formulation on \eqref{eq:relaxed-pde} and we obtain that $u$ and $\mathcal Z$ satisfy the system of equations in \eqref{eq:rof-pde}.

Finally, as in \cite[Theorem 2]{GiacomelliLasicaMoll2019} we note that, for $\varphi\in C^1(\overline{\Sigma})$,
	\begin{align*}
		&\int_{\Sigma}{\rm div}_g\big(\varphi(\mathcal B^\alpha+\sigma\nabla u^\alpha)\big)\,d\mu_g=\int_{\partial\Sigma}\varphi \<\mathcal B^\alpha+\sigma\nabla u^\alpha, \nu\>\,d\mu_{\tilde g},
\end{align*}
where $\tilde g$ is the induced Riemannian metric on $\partial \Sigma$, and
\begin{align*}
		&\int_{\Sigma}{\rm div}_g(\varphi(a_{\varepsilon_k}^\alpha))\,d\mu_g =0.
	\end{align*}
	Consequently, we have that
	\[
	\left|\int_{\partial\Sigma}\varphi \<\mathcal B^\alpha+\sigma\nabla u^\alpha, \nu\>\,d\mu_{\tilde g}\right|=\left|\int_\Sigma (h - h_k)\,d\mu_g\right|
	\]
with $h={\rm div}(\varphi(\mathcal B^\alpha+\sigma\nabla u^\alpha))$ and $h_k={\rm div}(\varphi(a_{\varepsilon_k}^\alpha))$. Then, since $h_k\to h$ in $L^2(\Sigma)$, we conclude that the boundary condition in \eqref{eq:rof-pde} holds.

\end{proof}

\begin{bem}
  Notice that with the same proof, but keeping $\epsilon>0$ fixed and letting $\delta,\sigma\to 0$, we also obtain Lipschitz regularity for minimizers of $\mathcal E_{0,\varepsilon}$. Additionally, the local statement in Corollary \ref{cor2} follows by redoing similar computations, but choosing $\varphi$ in Lemma \ref{Porr-Lem2} as a suitable cut-off function within a fixed ball.
\end{bem}

\subsection{Lipschitz regularity for the steady Mosolov problem} \label{LipMos}

In this part the goal is to prove claim (b) in Theorem \ref{main-thm3}. As before, the trick is to multiply $w$ by a suitable test function, which now has to depend on the target manifold in order to compensate the quadratic terms in $w$ arising from the positive upper curvature bound. Hereafter we will use the following notations:
\[\Phi = \frac{\varphi}{\psi \circ \rho }, \quad \text{with} \quad \rho =r_h \circ u := d_h(p, u(\cdot)), \quad \text{and} \quad \widetilde w = \Phi w,\]
where $\varphi = e^{c \eta}$ is the test function coming from the proof of Proposition \ref{prop:lipschitz} and $\psi$ will be specified later on. Let us first compute the terms appearing in \eqref{Ltest} with this new test function. Unless otherwise stated, all Laplacians, gradients and inner products are computed with respect to $g$, and $\psi$ is always evaluated in $\rho$. Hereafter, $C_{_\Sigma}$ does not have a specific meaning as before, it just indicates any constant which only depends on the domain surface.

\begin{lemma} \label{lem:PC-2}
With the previous notation, let $\kappa > 0$ be an upper bound for the absolute value of the sectional curvatures of $B_h(p, R) \subset \mathcal N$. Set $w = |d u|^2$ for $u$ any solution of \eqref{relaxed-pde-intr}. Then $\widetilde w$ attains its maximum at an interior point of $\Sigma$. Moreover, at this point,  the elliptic operator $\mathscr L$ defined as in \eqref{def-L} for $w > 1$ satisfies
\begin{align*}
\psi \mathscr L \widetilde w   \leq  &  \varphi \, \mathscr Lw +  \widetilde{w} \bigg(\frac1{\sqrt{w + \varepsilon^2}}  + \sigma \bigg) \big(\psi' \text{\rm co}_\kappa(\rho) w + \psi''|\nabla \rho|^2\big) - 2 \psi' \lambda R \widetilde w + \sigma \varphi \big|\nabla^{g\boxtimes h} d u\big|^2
\\ & + C_{_\Sigma}  \bigg(\frac{w^{3/2}}{\psi^2} \big(1- \psi' + (\psi')^2 + |\psi''|\big) - \psi' {\rm co}_{-\kappa}(\rho) \sqrt{w} + \frac1{4 \sigma}  \frac{(\psi')^2}{\psi}\bigg),
\end{align*}
provided that $0 <\psi \leq 1$ is a non-increasing function.
\end{lemma}

\begin{proof}
After applying the chain rule several times and some algebraic manipulation, we get
\begin{align*}
\Delta \Phi = \frac{\Delta \varphi}{\psi} - \frac{\Phi}{\psi}\big(\psi' \Delta \rho + \psi'' |\nabla \rho|^2\big) - 2 \frac{\psi'}{\psi} \<\nabla \rho, \nabla \Phi\>.
\end{align*}
On the other hand, let us check that \eqref{bdry-ineq2} still holds in the current setting. In fact, choosing again $x^2 = r_{_\partial}$ so that $\nabla_\nu = -\partial_2$, we obtain
\[\frac{\partial \widetilde w}{\partial \nu} = \frac{\partial}{\partial \nu}\bigg(\frac{e^{c \eta} w}{\psi \circ \rho}\bigg) = \varphi \, w \Big(- k_{_\partial} - c + \frac{\psi'}{\psi^2} \partial_2 \rho \Big) \leq  \widetilde w \, \frac{|\psi'|}{\psi} \frac{\partial{r_h}}{\partial y^\alpha} \partial_2 u^\alpha = 0,\]
where we have applied the Neumann boundary condition and the same choice of $c$ from Proposition \ref{prop:lipschitz}. Once more by virtue of Hopf's maximum principle, we can argue from now on at a maximum interior point for $\widetilde w$, at which as in \eqref{grad-max} we have
\begin{equation} \label{grad-w}
\nabla w = - w \nabla(\log \Phi).
\end{equation}

Next, using that
\begin{equation} \label{grad-Phi}
\nabla \Phi = \frac{\nabla \varphi}{\psi} - \frac{\Phi}{\psi} \psi' \nabla \rho,
\end{equation}
leads to
\[2\<\nabla w, \nabla \Phi\> = -2 w \big<\frac{\nabla \varphi}{\psi}, \nabla(\log \Phi)\big> +  2 w \frac{\psi'}{\psi} \<\nabla \rho, \nabla \Phi\>.\]
Summing up, we have
\begin{align*}
- \psi \Big(\Delta \Phi + \frac{2}{w}\<\nabla w, \nabla \Phi\>\Big) = 2 \frac{|\nabla \varphi|^2}{\varphi} - \Delta \varphi + \Phi \big(\psi' \Delta \rho + \psi'' |\nabla \rho|^2\big) - 2 \frac{\psi'}{\psi} \<\nabla \varphi, \nabla \rho\>.
\end{align*}

Doing similar computations for the remaining terms in  \eqref{Ltest}, we end up with
\begin{align*}
\mathscr L \widetilde w   & =   \Phi \, \mathscr Lw +  \frac{w}{\psi}\bigg(\frac1{\sqrt{w + \varepsilon^2}}  + \sigma \bigg) \bigg(2 \frac{|\nabla \varphi|^2}{\varphi} - \Delta \varphi + \Phi \big(\psi' \Delta \rho + \psi'' |\nabla \rho|^2\big) - 2 \frac{\psi'}{\psi} \<\nabla \varphi, \nabla \rho\>\bigg) \\ &  + \frac{w/\psi}{(w+ \varepsilon^2)^{\frac3{2}}} \bigg(\!\big(\nabla^2 \varphi -\Phi \psi' \nabla^2 \rho)(\nabla u^\alpha, \nabla u^\alpha) - \Phi \psi''\<\nabla \rho, \nabla u^\alpha\>^2 - 2\big<\nabla \varphi, \nabla u^\alpha\big>\big<\nabla u^\alpha, \frac{\nabla \Phi}{\Phi}\big>\bigg),
\end{align*}
where we have taken normal coordinates around $u(x) \in \mathcal N$, being $x$ the maximum interior point of $w$, and we understand sum over repeated indices.

By choosing $\psi$ to be a positive non-increasing function, using \eqref{bound-test} and taking into account that $|\nabla \rho|\leq \sqrt{w}$, we can estimate
\begin{align*}
\mathscr L \widetilde w   \leq  &  \Phi \, \mathscr Lw +  \frac{w}{\psi}\bigg(\frac1{\sqrt{w + \varepsilon^2}}  + \sigma \bigg) \bigg( C_{_\Sigma} \varphi + \Phi \big(\psi' \Delta \rho + \psi'' |\nabla \rho|^2\big) - 2 \frac{\psi'}{\psi} |\nabla \varphi|\sqrt{w}\bigg) \\ & -\frac{w/\psi}{(w+ \varepsilon^2)^{\frac3{2}}} \Phi \psi' \nabla^2 \rho(\nabla u^\alpha, \nabla u^\alpha)  + \frac{1/\psi}{(w+ \varepsilon^2)^{\frac1{2}}} \bigg(\! C_{_\Sigma} \varphi \, w + \Phi |\psi''|w^2 + 2 \psi \frac{|\nabla \varphi|}{\varphi} |\nabla \Phi| w\bigg)
\end{align*}

On the other hand, by the composition rule and the Hessian comparison theorem as in \cite[Lemma 2.9]{Sakai} for sectional curvatures bounded above by $\kappa$, we get
\begin{align*}
\psi' \Delta \rho = \psi' \Delta(r_h \circ u) =\psi' \Big(\nabla^2 r_h(\nabla u^\alpha, \nabla u^\alpha) +\big<^h\nabla r_h|_u, \tau(u)\big>_h\Big) \leq \psi' \text{co}_\kappa(\rho) w - \psi' |\tau(u)|_h.
\end{align*}
In turn, our equation \eqref{eqn-intr-2} allows us to get a bound for the last term above as follows
\begin{align*}
\bigg(\frac1{\sqrt{w + \varepsilon^2}}  + \sigma \bigg) \widetilde w \, |\tau(u)|_h & \leq  \widetilde w \bigg(\frac1{2(w+ \varepsilon^2)^{\frac3{2}}} |\nabla w| \sqrt{w}  +  \lambda \big|\exp_u^{-1} f|_h \bigg) \leq  \frac{\Phi}{2} |\nabla w|  + 2 R \lambda \widetilde w.
\end{align*}

It remains to estimate the term with $\nabla^2 \rho$, for which we use that the sectional curvatures of $\mathcal N$ are bounded below within the geodesic ball $B_h(p, R) \supset u(\overline \Sigma)$, say by  $- \kappa$. Then the remaining bound of the Hessian comparison theorem leads to
\begin{align*}
\frac{|\psi'|}{\psi} \Phi \nabla^2 \rho(\nabla u^\alpha, \nabla u^\alpha) &\leq \frac{|\psi'|}{\psi} \Phi \Big(\text{co}_{-\kappa}(\rho) w + \big|\nabla^{g\boxtimes h} d u\big|\Big) \\ & \leq  \frac{|\psi'|}{\psi}  \text{co}_{-\kappa}(\rho) \widetilde{w} + \sigma \Phi \big|\nabla^{g\boxtimes h} d u\big|^2 + \frac1{4 \sigma} \Phi \bigg(\frac{\psi'}{\psi}\bigg)^2,
\end{align*}
which follows by Young's inequality.

The last ingredient to bound $\mathscr L \widetilde w$ at an interior maximum is that
\begin{align} \label{nablaw}
|\nabla w| = \frac{w}{\Phi} |\nabla \Phi| \leq w\bigg(\frac{|\nabla \varphi|}{\psi \Phi} +  \frac{|\psi'|}{\psi} |\nabla \rho|\bigg) \leq w \frac{|\nabla \varphi|}{\varphi} -\frac{\psi'}{\psi} w^{3/2},
\end{align}
where we have applied \eqref{grad-w} and \eqref{grad-Phi}.

Gathering all the above inequalities, taking into account that $\varphi$ is by construction a bounded function on $\Sigma$, and that we can assume $w > 1$ without loss of generality   (as otherwise $w$ is bounded and we are done), after a lengthy but straightforward manipulation we reach the claimed inequality.

\end{proof}

\begin{lemma} \label{Lw-PCC}
	With the previous notation, let $\kappa > 0$ be an upper bound for the absolute value of the sectional curvatures of $B_h(p, R) \subset \mathcal N$ with $R < R_\kappa^\ast$. Set $w = |d u|^2$ for $u$ any solution of \eqref{relaxed-pde-intr}. For the elliptic operator $\mathscr L$ defined as in \eqref{def-L} one has
	\begin{align*}
	\mathscr L w   + 2 \sigma   \big|\nabla^{g\boxtimes h} d u\big|^2_{g\boxtimes h}    \leq &  \bigg(\!\frac1{\sqrt{w + \varepsilon^2}}  + \sigma\! \bigg)   (\kappa \, w +   C_{_\Sigma}) w    +  2\lambda   \pi |df_\delta| \sqrt{w} \\ \nonumber \smallskip & +  \frac{|\nabla w|^2_g}{2(w +\varepsilon^2)^{\frac3{2}}} + 2 \lambda R \frac{|\nabla w|_g}{\sqrt{w + \varepsilon^2}}.
	\end{align*}
\end{lemma}

\begin{proof}
From the proof of Lemma \ref{eqn-intr-2}, as we have only used the curvature assumption to estimate the term with $\lambda$, we have the equality
\begin{align*}
\mathscr L w   + 2 \bigg(\frac1{\sqrt{w + \varepsilon^2}}  + \sigma \bigg) & \bigg(\big|\nabla^{g\boxtimes h} d u\big|^2_{g\boxtimes h}  - \mathscr R  \bigg)  - 2 \lambda  \big<\widetilde \nabla  \big(\exp^{-1}_u f_\delta\big), du\big>_{g\boxtimes h} \\ & =  \frac3{2} \frac{\big<dw \otimes \<\nabla w, \nabla u\>_g, du\big>_{g\boxtimes h}}{(w +\varepsilon^2)^{5/2}}   -  \frac{\<dw \otimes \tau(u), du\>_{g\boxtimes h}}{(w +\varepsilon^2)^{3/2}}   - \frac{|\nabla w|^2_g}{2(w +\varepsilon^2)^{3/2}},
\end{align*}

A slight variation of the argument from the proof of Lemma \ref{eqn-intr-2} which uses Jacobi field comparison, so that it works for positive curvature bounds (after reparametrization to get a unit speed geodesic and some slight modification for non necessarily normal fields) and by means of Hessian comparison for the squared distance function (see \cite[Theorem 6.6.1]{Jost2017}), we get 	
\begin{align}\label{jacobi-comp-general}
-  \big<\widetilde \nabla  \big(\exp^{-1}_u f_\delta\big), du\big>_{g\boxtimes h} \geq d_h(u, f) \, {\rm co}_\kappa(d_h(u, f)) \, w \, - 2 \frac{d_h(u, f)}{{\rm s}_\kappa(d_h(u, f))} |df_\delta| \sqrt{w} \geq - \pi |df_\delta| \sqrt{w},
\end{align}
where to discard the first term we have applied that $d_h(u,f) < 2 R < \frac{\pi}{2\sqrt{\kappa}}$ by the definition of $R_\kappa^\ast$ in \eqref{Rast}.

Combining the latter with \eqref{term52}, we can write
\begin{align} \label{aux-Lw-curvP}
\mathscr L w   + 2 \bigg(\frac1{\sqrt{w + \varepsilon^2}}  + \sigma \bigg)  \big|\nabla^{g\boxtimes h} d u\big|^2_{g\boxtimes h}    & \leq 2 \bigg(\frac1{\sqrt{w + \varepsilon^2}}  + \sigma \bigg)   \mathscr R    + 2 \lambda   \pi |df_\delta| \sqrt{w} \nonumber \smallskip \\ &   -  \frac{\<dw \otimes \tau(u), du\>_{g\boxtimes h}}{(w +\varepsilon^2)^{3/2}}   + \frac{|\nabla w|^2_g}{2(w +\varepsilon^2)^{3/2}}.
\end{align}

For the curvature term given by \eqref{curv-BW}, we have
\begin{align*}
\mathscr R \leq \kappa \sum_{i,j }\big|u_\ast e_i \wedge u_\ast e_j\big|_h^2 + C_{_\Sigma} w \leq \kappa \sum_{i,j }\Big(|u_\ast e_i|^2_h  |u_\ast e_j|^2_h - \<u_\ast e_i, u_\ast e_j\>_h^2\Big) + C_{_\Sigma} w \leq \frac{\kappa}{2} w^2 + C_{_\Sigma} w,
\end{align*}
which follows thanks to the standard estimates
\[\sum_{i,j}\<u_\ast e_i, u_\ast e_j\>_h^2 \geq \sum_i |u_\ast e_i|^4_h \geq \frac1{2} \bigg(\sum_i |u_\ast e_i|^2_h\bigg)^2 = \frac1{2} w^2.\]

In turn, by applying Kato's inequality, which yields $|\nabla w|_g \leq 2 \sqrt{w} \big|\nabla^{g\boxtimes h} d u\big|_{g\boxtimes h}$,  to the last term in \eqref{aux-Lw-curvP}, we obtain
\begin{align*}
\mathscr L w   + 2 \sigma   \big|\nabla^{g\boxtimes h} d u\big|^2_{g\boxtimes h}    & \leq  \bigg(\frac1{\sqrt{w + \varepsilon^2}}  + \sigma \bigg)   (w^2 + 2\,  C_{_\Sigma} w)    + \lambda   \pi |df_\delta| \sqrt{w} \nonumber \smallskip +  \frac{|\nabla w|_g\, |\tau(u)|_h}{w +\varepsilon^2}.
\end{align*}

Finally, our equation \eqref{eqn-intr-2} leads to
\begin{align*}
|\tau(u)|_h \leq \bigg(\frac1{\sqrt{w + \varepsilon^2}}  + \sigma \bigg)^{-1}\bigg(\frac{|\nabla w|_g\sqrt{w}}{2(w + \varepsilon^2)^{3/2}} + \lambda \, d_h(u,f)\bigg) \leq \frac{|\nabla w|_g}{2(w + \varepsilon^2)^{1/2}} + 2 \lambda R \sqrt{w + \varepsilon^2},
\end{align*}
which leads to the inequality in the statement.
\end{proof}

By combining the estimates in the previous two lemmas, after some extra work and an appropriate choice of the test function $\psi$, we can prove a Lipschitz regularity result for constants uniform in $\varepsilon$, but depending on $\sigma$.

\begin{prop}\label{prop:lipschitz-posit}
	Let $u$ be a minimizer of $\mathcal E_{\sigma, \varepsilon}$ when $\kappa > 0$ is an upper bound for the absolute value of the sectional curvatures in $B_h(p,R) \subset \mathcal N$ for $R < R_\kappa^\ast$. Then there exist constants $C >0$, depending on the domain surface $\Sigma$, $\kappa$, $\lambda$ and/or $R$ (but are independent of $\varepsilon, \sigma$), such that the following bound
	\[
	\sigma^2 |d u| \leq  C \big({\rm Lip} (f) + 1\big)
	\]
	holds.
\end{prop}

\begin{proof}
	Setting again $\widetilde w = \Phi w$ and with the notation at the beginning of this subsection, Lemma  \ref{lem:PC-2} ensures that $\widetilde w$ attains its maximum at an interior point $x$. We perform all the computations at $x$ hereafter. Let us multiply by $\Phi$ the inequality from Lemma \ref{Lw-PCC}, which joint with the one in Lemma \ref{lem:PC-2} and further estimates using \eqref{nablaw}, allow us to write
	\begin{align*}
0 \leq \mathscr L \widetilde w   +  \Phi \sigma   \big|\nabla^{g\boxtimes h} d u\big|^2_{g\boxtimes h}     \leq & \frac{\Phi}{\psi} \bigg(\!\frac1{\sqrt{w + \varepsilon^2}}  + \sigma\! \bigg)   \Big((\psi' \text{\rm co}_\kappa(\rho) + \kappa \psi) w^2 +   (C_{_\Sigma} \psi + \psi'' |\nabla \rho|^2) w\Big)   \\   & +  2 \lambda \Phi  \pi |df_\delta| \sqrt{w}  - 4 \frac{\psi'}{\psi} \lambda R \widetilde w  - \frac{\psi'}{\psi} {\rm co}_{-\kappa}(\rho) C_{_\Sigma} \sqrt{w} + 2 \lambda R C_{_\Sigma} \Phi\sqrt{w}
\\ & +  \frac{C_{_\Sigma}}{\psi^2} w^{3/2} \bigg(1- \psi'  + (\psi')^2 + |\psi''|\bigg) + \frac1{4 \sigma} C_{_\Sigma} \frac{(\psi')^2}{\psi^2}.
\end{align*}

	As, by assumption,  $\rho(x) < \frac{\pi}{4 \sqrt{\kappa}}$, set $\psi = {\rm c}_\kappa^2$. Taking into account that $\frac{\psi'}{\psi} = -2\kappa \, {\rm ta}_\kappa$ and  $\psi'' \leq 0$, we can estimate	
\begin{align*}
\sigma \Phi \, \text{\rm{c}}_\kappa^2(\rho)  w^2 \leq C_{_{\lambda, \Sigma}} |df_\delta| \sqrt{w} + C (w^{3/2} + \sigma^{-1})
\end{align*}
where $C$ is a constant that may depend on $\Sigma$, $\lambda$, $\kappa$ and/or $R$, but is independent of $\varepsilon, \sigma$ and has a meaning that may change from line to line. Notice that we have used that
\[\frac{|\psi'|}{\psi} {\rm co}_{-\kappa}(\rho) \leq 2 \sin({\sqrt{\kappa} \rho})\bigg(1 + \frac1{\sqrt{\kappa}\rho}\bigg) \leq C.\]

Accordingly, we conclude
\begin{align*}
\sigma  |d u| \leq C_{_{\lambda, \Sigma}} |df_\delta|  + C (1 + \sigma^{-1}),
\end{align*}
which leads to the Lipschitz regularity claimed.
\end{proof}

The same proof as that of Proposition \ref{final_result}, together with lower semicontinuity of $\mathcal E_\sigma$ yields the following regularity result without further restrictions on the curvature of the manifold $\mathcal N$:
\begin{prop}
		For any $\sigma > 0$ there exists a minimizer $u$ of $\mathcal E_\sigma$, such that $u\in C^{0,1}(\overline\Sigma; \mathcal N)$.
	\end{prop}

\section{Regularity for signal denoising} \label{signal}

	In this section, we will prove Lipschitz regularity of minimizers of $\cE_{\varepsilon,\sigma}$, for $\varepsilon,\sigma\geq0$ and $f\in \text{Lip}(\Gamma)$. We will only detail the proofs for the case of $\Gamma=[0,1]$, the situation of $\Gamma=\mathbb{S}^1$ being simpler because of the absence of boundary conditions. In this setting, the system of Euler-Lagrange equations reads as (again we identify manifold-valued functions with the embedded representation given by Nash theorem).
\begin{equation}\label{eq:system-1d}
\left\{
\begin{array}{rcll}
(\cZ_{\varepsilon, \sigma}^\alpha(u'))' + \mathcal (\mathcal A_u(u', \cZ_{\varepsilon,  \sigma}(u')))^\alpha  &=&  -\lambda (\exp_u^{-1} f)^\alpha \quad &\text{in $(0,1)$}, \smallskip \\
(u')^\alpha &=& 0 \quad &\text{at $\{0,1\}$}
\end{array}\tag{$\mathcal S_{\varepsilon, \sigma}^{f}$},
\right.
\end{equation}
with $\cZ_{\varepsilon, \sigma}^\alpha(u')=\Big(\frac{1}{\sqrt{\varepsilon^2+|u'|^2}}+\sigma \Big)(u')^\alpha$.

In this case, we first show  that weak solutions have better a-priori regularity that in the case of the domain being a surface $\Sigma$.

\begin{lemma}\label{lem:1dreg}
	Let $u$ be a bounded weak solution to \eqref{eq:system-1d} for $\varepsilon,\sigma>0$. Then $u'\in C^{0,1}([0,1])$.
\end{lemma}

\begin{proof}
 First of all, we observe that since $u\in H^1(0,1)$, direct from the system, one obtains $\cZ_{\varepsilon,\sigma}(u')\in W^{1,1}(0,1)$.

	We consider now the transformation $\Phi:\R^N\to\R^N$ given by $$\Phi(\xi)=\Big(\frac{1}{\sqrt{\varepsilon^2+|\xi|^2}}+\sigma\Big)\xi.$$ Observe that  $\Phi^{-1}$ is well-defined and belongs to $C^1(\R^N)$. Therefore, as $u' = \Phi^{-1}(\cZ_{\varepsilon,\sigma}(u'))$ we obtain that $u'\in W^{1,1}(0,1)$. Hence, by Sobolev embedding, $u'$ is bounded. Then, once again from the system, we get that $\cZ_{\varepsilon,\sigma}(u')$ is a Lipschitz function, which implies that
 $u' $ is also Lipschitz.
\end{proof}

We are now in position to prove Theorem \ref{thm:1dcase}. 

\begin{proof}[Proof of Theorem \ref{thm:1dcase}] First of all, existence of a minimizer follows easily as in the proof of Proposition \ref{prop:existence}, by using Remark \ref{rem:1dcase}. As the proof of  item $(a)$ mimics that of Theorem \ref{main-thm1}, we only give a brief sketch  since all the arguments have been already shown in a more complicated framework.

  Once the regularity of the solutions of the approximating problems has been obtained in Lemma \ref{lem:1dreg}, a bootstrap argument and mollification of the data $f$ as in Proposition \ref{prop:holder-D2u} yield smoothness of the approximating minimizers. Then, the same Bernstein technique as in Proposition \ref{Porr-cvx}, shows that indeed the Lipschitz constants of the approximations are uniformly bounded. We point out that no curvature restrictions appear in the 1-dimensional case, since the term $\mathscr R$ in Lemma \ref{Porr:lem1} vanishes. Finally, a passage to the limit as in Proposition \ref{final_result}, together with lower semicontinuity of $\mathcal E_\sigma$, proves the result.

  To prove item $(b)$ we proceed as in \cite{GL19} (see also \cite{GiacomelliLasicaMoll2019, GroLa}). We let $|\cdot|_\varepsilon:=\sqrt{\varepsilon^2+|\cdot|^2}$. First, we will show that for every $\tilde\Gamma \Subset \Gamma$, it holds
		\begin{equation}\label{eq:lp-compactly}
			\int_{\tilde\Gamma} |u'|^p_\varepsilon d x \leq {C^p} \int_{\tilde \Gamma}|f'|^p d x  + p \,\varepsilon^p|\tilde \Gamma| \quad \text{for $p\in[1,2]$},
		\end{equation}
		where $C > 0$, and
		\begin{equation}\label{eq:lp-ball}
			\limsup_{\varepsilon \to 0} \int_{B(x_0,r)} |u'|_\varepsilon^p d x  \leq {C^p} \int_{B(x_0, R)} |f'|^p d x  \quad \text{for $p\in(1,2]$},
		\end{equation}
		for $0 < r < R$ and $x_0\in \tilde \Gamma$  such that $B(x_0,R) \Subset \Gamma$.

We let $u$ be a solution to system $(\mathcal S^f_{\varepsilon, \varepsilon^2})$, which can be rewritten as $$\nabla_t\cZ_{\varepsilon, \sigma}^\alpha(u') =- \lambda \exp_u^{-1} f,$$ with $\nabla_t$ denotes the covariant derivative along the curve $t\mapsto u(t)$. By the regularity obtained in Lemma \ref{lem:1dreg}, we can differenciate the system and given $\varphi$ any Lipschitz function, we take
		$
		\bar\varphi := -\varphi^2 |u'|_\varepsilon^{p-2} u'$,
		as a test function. Therefore,
		\begin{equation}\label{eq:main-eq}
			\lambda \int_{\tilde \Gamma} \langle\bar\varphi, (\exp_u^{-1}(f))'\rangle d t  = \int_{\tilde \Gamma} \langle \bar\varphi',\nabla_t {\mathcal Z_{\varepsilon,\varepsilon^2}(u')}\rangle d t .
		\end{equation}
		The right-hand side of the above equation is equal to
		\begin{align}
			-\int_{\tilde \Gamma} (\varphi^2 |u'|_\varepsilon^{p-1})' \big\langle \frac{u'}{|u'|_\varepsilon}, \nabla_t {\mathcal Z_{\varepsilon,\varepsilon^2}(u')}\big\rangle  d t  -\int_{\tilde \Gamma} \varphi^2 |u'|_\varepsilon^{p-1} \big\langle\left(\frac{u'}{|u'|_\varepsilon}\right)', \nabla_t {\mathcal Z_{\varepsilon,\varepsilon^2}(u')}\big\rangle  d t =I_1+I_2 . \label{eq:our-terms}
		\end{align}
		
		\bigskip
		Taking into account that \mbox{$\langle u', \nabla_t {\mathcal Z_{\varepsilon,\varepsilon^2}(u')}\rangle= \langle u', \big(\mathcal Z_{\varepsilon,\varepsilon^2}(u')\big)'\rangle$}, we obtain $$I_1=-2\varepsilon^2\int_{\tilde\Gamma}\varphi\varphi'|u'|^{p-5}_\varepsilon(1+|u'|^3_\varepsilon)\langle u',u''\rangle dt-(p-1)\varepsilon^2\int_{\tilde\Gamma}\varphi^2|u'|^{p-7}_\varepsilon(1+|u'|^3_\varepsilon)\langle u',u''\rangle^2 dt.$$ 
On the other hand, for the inner product in $I_2$ we have
\begin{align*} \big\langle\left(\frac{u'}{|u'|_\varepsilon}\right)', \nabla_t {\mathcal Z_{\varepsilon,\varepsilon^2}(u')}\big\rangle  & = \big\langle\nabla_t\frac{u'}{|u'|_\varepsilon}, \nabla_t\frac{u'}{|u'|_\varepsilon} + \varepsilon^2 \nabla_t u'\big\rangle \\ & = \bigg|\nabla_t\frac{u'}{|u'|_\varepsilon}\bigg|^2 + \frac{\varepsilon^2}{|u'|_\varepsilon}\bigg(\big|\nabla_t u'\big|^2 - \frac{\langle\nabla_t u', u'\rangle^2}{|u'|_\varepsilon^2}\bigg)
 \geq 0. \end{align*}

Thus, \begin{equation}
  \label{aux:gl} \frac{\lambda}{\varepsilon^2} \int_{\tilde \Gamma} \langle\bar\varphi, (\exp_u^{-1}(f))'\rangle d t\leq -\int_{\tilde\Gamma}\varphi|u'|^{p-7}_\varepsilon(1+|u'|^3_\varepsilon)\langle u',u''\rangle(2\varphi'|u'|^2_\varepsilon+(p-1)\varphi \langle u',u''\rangle) dt
\end{equation}

		Therefore,  letting $\varphi\equiv 1$, we deduce from \eqref{aux:gl} and \eqref{jacobi-comp-general} from Lemma \ref{Lw-PCC} that
		\begin{equation}\label{eq:inequality-u-f}
			\pi\int_{\tilde\Gamma} |f'||u'||u'|_\varepsilon^{p-2} dt- \tilde C \int_{\tilde\Gamma} |u'|^2|u'|_\varepsilon^{p-2} dt\geq \int_{\tilde \Gamma} |u'|^{p-2}_\varepsilon\langle u',(\exp_u^{-1}(f))'\rangle  d t\geq 0,
		\end{equation}
where $\tilde C:=\min\{2R {\rm co}_{\kappa}(2R),1\}$.
	
		Now, noting that for $p \leq 2$ we have $|u'|_\varepsilon^{p-2}|u'|^2 \geq |u'|_\varepsilon^p - \varepsilon^p$, we obtain
		\begin{equation*}\label{eq:u-part}
			\begin{aligned}
				\int_{\tilde \Gamma} |u'|_\varepsilon^p d t - \varepsilon^p|\tilde \Gamma| \leq \int_{\tilde \Gamma} |u'|_\varepsilon^{p-2} |u'|^2 d t
				\leq \,  	\frac{\pi}{\tilde  C}\int_{\tilde\Gamma} |f'||u'|_\varepsilon^{p-1} dt 
				\leq \frac{p-1}{p} \int_{\tilde \Gamma} |u'|_\varepsilon^p d t + \frac{\pi^p}{p \tilde  C^p}\int_{\tilde \Gamma} |f'|^p d t,
			\end{aligned}
		\end{equation*}
which follows using Young's inequality in the last expression. From here we reach \eqref{eq:lp-compactly}, and an argument identical to that of \cite[Lemma]{GL19} yields \eqref{eq:lp-ball}.

 Finally, by replicating the proof of \cite[Theorem]{GL19}, we obtain (b).
\end{proof}

\end{document}